\documentclass[a4paper,12pt,oneside,reqno]{article}
\usepackage{amsfonts,amssymb,amsmath,amsthm}
\usepackage[a4paper]{typearea}
\usepackage[utf8]{inputenc}

\newtheorem{theorem}{Theorem}[section]
\newtheorem{lemma}[theorem]{Lemma}
\newtheorem{proposition}[theorem]{Proposition}
\newtheorem{corollary}[theorem]{Corollary}
\newtheorem*{claim}{Claim}

\theoremstyle{definition}
\newtheorem{definition}[theorem]{Definition}
\newtheorem{assumption}[theorem]{Assumption}
\newtheorem{condition}{Condition}

\theoremstyle{remark}
\newtheorem*{remark}{Remark}

\numberwithin{equation}{section}
\long\def\rem#1{}

\def\dist{\mathop{\rm dist}\nolimits}

\def\bbone{{\mathchoice {\rm 1\mskip-4mu l} {\rm 1\mskip-4mu l}
          {\rm 1\mskip-4.5mu l} {\rm 1\mskip-5mu l}}}

\def\vagconv{\overset v\to }
\def\ppconv{\overset {pp} \to }

\def\laweq{\overset{\mathrm{law}}=}

\def\Asl{\mathop{\mathsf{Asl}}\nolimits}
\def\BTM{\mathop{\mathsf{BTM}}\nolimits}

\def\Id{\mathrm{Id}}
\def\d{\mathrm{d}}
\def\dd{\mathtt{d}}

\def\<{\langle}
\def\>{\rangle}

\def\mm#1#2{{#2}\circ{#1}}
\def\inp{a}

\begin{document}
\title{Dynamics of trap models}
\author{G\'erard Ben Arous \and Ji\v r\'{\i} \v Cern\'y}
%
\maketitle
\ \ 

\vskip -2.5cm

\ \ 
\tableofcontents

\section{Introduction} 

These notes cover one of the topics of the class given in the Les 
Houches Summer School ``Mathematical statistical physics'' in July 
2005. The lectures tried to give a summary of the 
recent mathematical results about the long-time behaviour of dynamics 
of (mean-field) spin-glasses and other disordered media (expanding on 
  the review \cite{Ben02}). We have 
chosen here to restrict the scope of these notes to the dynamics of 
trap models only, but to cover this topic in somewhat more depth.

Let us begin by setting the stage of this long review about the trap 
models by going back to one of the motivations behind their 
introduction by Bouchaud \cite{Bou92}, i.e.~dynamics of spin-glasses,
which is
 indeed a very good field to get an accurate sample of 
possible and generic long-time phenomena (as aging, memory, 
  rejuvenation, failure of the fluctuation-dissipation theorem;
  see~\cite{BCKM98} for a global review).

\smallskip

This class of problems can be roughly described as follows. Let 
$\Gamma$ (a compact metric space) be the state space for spins and 
$\nu$ be a probability measure on $\Gamma$. Typically, in the 
discrete (or Ising) spins context $\Gamma=\{-1,1\}$ and 
$\nu=(\delta_{1}+\delta_{-1})/2$. In the continuous or soft spin 
context $\Gamma=I$, a compact interval  of the real line, and 
$\nu(dx) = Z^{-1}e^{-U(x)}dx$, where $U(x)$ is the ``one-body 
potential''. For each configuration of the spin system, i.e.~for each 
$\boldsymbol x=(x_1,...,x_n) \in \Gamma^n$ one defines a random 
Hamiltonian, $H^n_{\boldsymbol J}(\boldsymbol x)$, as a function of 
the configuration $\boldsymbol x$ and of an exterior source of 
randomness $\boldsymbol J$, i.e.~a random variable defined on another 
probability space. The Gibbs measure at inverse temperature $\beta$ 
is then defined on the configuration space $\Gamma^n$ by
\begin{equation}
  \mu^n_{\beta,\boldsymbol J}(\d \boldsymbol x)
  = \frac 1{Z^n_{\boldsymbol J}}
  \exp\big(-\beta H^n_{\boldsymbol J}(\boldsymbol x)\big)
  \nu(\d\boldsymbol x).
\end{equation}
The statics problem amounts to understanding the large $n$ behaviour 
of these measures for various classes of random Hamiltonians 
(\cite{Tal03} is a recent and beautiful book on the mathematical 
  results pertaining to these equilibrium problems). The dynamics 
question consists of understanding the behaviour of Markovian 
processes on the configuration space $\Gamma^n$, for which the Gibbs 
measure is invariant and even reversible, in the limit of large 
systems (large $n$) and long times, either when the randomness 
$\boldsymbol J$ is fixed (the quenched case) or when it is averaged 
(often called the annealed case in the mathematics literature, but 
  not in the physics papers). These dynamics are typically Glauber 
dynamics for the discrete spin setting, or Langevin dynamics for 
continuous spins. 

Defining precisely what we mean here by large system size $n$ and 
long time $t$ is a very important question, and very different 
results can be expected for various time scales $t(n)$ as functions 
of the size of the system. The very wide range of time (and energy or 
  space) scales present in the dynamics of these random media is the 
main interest and difficulty of these questions (in our view). At one 
end of the spectrum of time scales one could take first the limit 
when $n$ goes to infinity and then $t$ to infinity. This is the 
shortest possible long-time scale, much too short typically to allow 
any escape from metastable states since the energy barriers the 
system can cross are not allowed to diverge. This short time scale is 
well understood for dynamics of various spin-glass models (and 
  related models) in the physics literature, in particular for the 
paradigmatic Langevin dynamics of spherical $p$-spin models of spin-glasses, 
mainly 
through the equations derived by Cugliandolo and Kurchan \cite{CK93} 
(see also \cite{CHS93}). 
The fact that the results given by this short-time limit are seen as 
correct for models with a continuous replica symmetry breaking (like 
  the Sherrington-Kirkpatrick model) is one of the many bewildering 
predictions made by the physicists. In the Les Houches lectures we 
covered some of the recent mathematical results about this short-time 
scale for Langevin dynamics obtained in collaboration with Alice 
Guionnet and  Amir Dembo (\cite{BDG06,BDG01,BG97,Gui97}). For lack of space and in order 
to keep a better focus we will not touch this topic here at all.

On the contrary we will be interested in the other end of the range 
of time scales, i.e.~time scales $t_w(n)$ depending on the system size 
in such a way that they allow for the escape from the deep metastable 
states. This is where the introduction of the phenomenological trap 
models by Bouchaud becomes meaningful. We will now try to explain the 
relevance of these models in this setting, although in necessarily 
rather imprecise terms. 

At low temperature the Gibbs measure $\mu_{\beta ,\boldsymbol J}^n$ 
should be essentially carried by a small part of the configuration 
space, the ``deep valleys'' of the random landscape, i.e.~the regions 
of low energy, (the ``lumps'' of $p$-spin models for instance). So 
that the dynamics should spend most of the time in these regions, 
which thus becomes very sticky ``attractors'' or ``traps''. The trap 
models ignore the details of the dynamics inside these sticky 
regions. They only keep the statistics of the height of the barriers 
that the system must cross before leaving these regions, and 
therefore, the statistics of the trapping times (i.e.~the times needed 
  to escape them). Moreover, the trap models keep the structure of 
the possible routes from one of these traps to the others. The 
dynamics is, therefore, reduced to its caricature: it lives only on 
the graph whose vertices are all the relevant attractors and the 
edges are the pairs of communicating attractors (see 
  Section~\ref{s:definition} for precise definition).

It is a non-trivial matter to prove that these phenomenological 
models could be of any relevance for the original problems. In fact, 
in the first introduction \cite{Bou92}of this model a large complete 
graph  was supposed to be a good ansatz for the simplest model of a 
mean-field spin-glass, i.e.~Derrida's Random Energy Model 
\cite{Der81}. Proving rigorously that Bouchaud's ansatz or 
``phenomenological'' model was indeed a very good approximation for 
metastability and aging questions is quite delicate, and was done 
only recently, initially in \cite{BBG03,BBG03b}. Later we realised 
(\cite{BC06}) that this very simple Bouchaud ansatz for the REM was 
in fact even better, in the sense that the range of time scales where 
it is a reliable approximation is very wide.

We now believe that this relevance is even much wider (``universal'', 
  if we dare) in the following sense: some of the lessons learnt from 
Bouchaud's picture in the REM should be relevant for very wide 
classes of mean-field spin-glasses in appropriate time scales. The 
art is in choosing these time scales long enough so that the (usually 
  sparse) deep traps can be found and thus some trapping can take 
place, but short enough so that the deepest traps are not yet 
relevant and thus the equilibrium (which is of course heavily 
  model-dependent) is not yet sampled by the dynamics. This belief 
has propped us into understanding more deeply the trap models in the 
larger possible generality of the graph structure and of the time 
scales involved. We have indeed found a very universal picture valid 
for all the examples we have studied, except for the very particular 
one-dimensional Bouchaud trap model which belongs to another class, 
as we will see below.

\smallskip

We want to emphasise here that the dynamics of the mean-field 
spin-glasses is far from being the only motivation that makes the 
study of trap models worthwhile, see e.g.~\cite{BB03} for references 
of applications to fragile glasses, soft glassy and granular 
materials, and pinning of extended defects. 

\smallskip

Let us now describe what these notes contain in more detail. We start 
by giving in Section \ref{s:definition} the definition of the 
Bouchaud trap model for a general graph and a general ``depth'' 
landscape. We then study, in Section \ref{s:onedim}, the very 
specific one-dimensional case (the graph here is $\mathbb Z$) in much 
detail. The results and the methods are  different from all other 
situations we study. We rely essentially on the scaling limit 
introduced by Fontes, Isopi and Newman in \cite{FIN02}. This scaling 
limit is an  interesting self-similar singular diffusion, which  
gives quite easily results about aging, subaging and the 
``environment seen from the particle''. We then go, in 
Section~\ref{s:ddim}, to the $d$-dimensional case and show that the 
essence of the results is pretty insensitive to the dimension (as 
  long as $d\ge 2$, with some important subtleties for the most 
  difficult case, i.e.~$d=2$). In particular we also give a  scaling 
limit, quite different from the one-dimensional case. We show that 
the properly-rescaled ``internal clock'' of the dynamics converges to an 
$\alpha$-stable subordinator and that the process itself when 
properly rescaled converges to a ``fractional-kinetics'' type of 
dynamics (\cite{Zas02}) which is simply the time change of a $d$-dimensional 
Brownian Motion by an independent process, the inverse of an $\alpha$-stable 
subordinator. This process is also a self-similar continuous process, 
but it is no longer Markovian. In fact, this scaling-limit result is in some 
sense a (non-trivial) triviality result. It says that the Bouchaud 
trap model has the same scaling limit as a Continuous Time Random 
Walk \`a la Montroll-Weiss \cite{MW65}. The aging results are then 
seen as a direct 
consequence of the generalised arcsine law for stable subordinators.

This picture (valid for all $d\ge 2$) is also naturally valid for the 
infinite-dimensional (or mean-field) case, i.e.~for large complete 
graphs  which we study in Section~\ref{s:general}. Thus, we see that 
the critical mean-field dimension is $2$ (in fact, we do not really 
  guess what could happen for dimensions between $1$ and $2$, but it 
  could be an interesting project to look at these models on say 
  deterministic fractals with spectral dimension between $1$ and $2$). 
For large complete graphs, which is a very easy case, we choose to 
give a new proof instead of following the well-established route 
using renewal arguments (as in  \cite{BD95} and in 
  \cite{BBG03b}).   This proof is slightly longer but illustrates, in this 
simple context, the strategy that we follow in other more difficult 
cases. An advantage of  this line of proof  is worth mentioning:  
we  get aging results in longer time scales 
than usual.

We then use the intuition we hope to have given in the easy Section 
\ref{ss:complete} to explain in Section~\ref{ss:universal} a general 
(universal?) scheme for aging based on the same arguments. We isolate 
the arguments needed for the proof given for the complete graph to 
work in general. This boils down to six technical conditions under 
very general circumstances. We then show how, under usual 
circumstances, these conditions can be reduced to basic 
potential-theoretic conditions for the standard random walk and 
random subsets of the graph. This general scheme is shown to be 
applicable not only to the cases we already know (i.e.~the Bouchaud 
  model on $\mathbb Z^d$ with $d\ge 2$, or the large complete 
  graphs), but also to dynamics of the Random Energy Model (in a wide 
  range of time scales, shorter than the one given in \cite{BBG03b}, 
  including some above the critical temperature) and also to very long 
time scales for large boxes in finite dimensions (with periodic 
  boundary conditions).


\section{Definition of the Bouchaud trap model}
\label{s:definition}
We define here a general class of reversible Markov chains on graphs, 
which were introduced by Bouchaud \cite{Bou92}  in order to 
give an effective model for trapping phenomena. The precise 
definition of these Markov chains necessitates three ingredients: a 
graph $G$, a trapping landscape $\boldsymbol \tau $ and an extra real 
parameter $a$. 

We start with the graph $G$,  $G=(\mathcal V, \mathcal E)$, 
with the set of vertices  $\mathcal V$ and   
the set of edges $\mathcal E$. We suppose that $G$ is non-oriented and 
connected; $G$ could be finite or infinite. 

We then introduce the trapping landscape $\boldsymbol \tau $. For 
each vertex $x$, $\tau_x$ is a positive real number which  is
referred to as the depth of the trap at $x$. $\boldsymbol \tau $ 
can also be seen as a (positive) measure on $\mathcal V$,
\begin{equation}
  \boldsymbol \tau=\sum_{x\in \mathcal V} \tau_x\delta_x.
\end{equation}

Finally, we define the continuous time Markov chain $X(t)$ on 
$\mathcal V$ by its jump rates $w_{xy}$,
\begin{equation}
  \label{e:rates}
  w_{xy}=\begin{cases}
    \nu \tau_x^{-(1-a)}\tau_y^a,&\text{if $\< x,y\> \in \mathcal E$,}\\
    0,&\text{otherwise.}
  \end{cases}
\end{equation}
Hence, the generator of the chain is 
\begin{equation}
  \label{e:generator}
  Lf(x)=\sum_{y:\<x,y\>\in \mathcal E}w_{xy}\big(f(y)-f(x)\big).
\end{equation}
Here, the linear scaling factor $\nu $ defines a time unit and is 
irrelevant for the dynamical properties. Its value in these notes 
varies for different graphs, mostly for technical convenience. The 
parameter $a\in[0,1]$ characterises the ``symmetry'' or ``locality'' 
of the dynamics; its role will be explained later. The initial state 
of the process will be also specified later, usually we set 
$X(0)=\boldsymbol 0$, where $\boldsymbol 0$ is an arbitrary fixed 
vertex of the graph.
\begin{definition}
  Given a  graph $G=(\mathcal V, \mathcal E)$, a trapping 
  landscape  $\boldsymbol \tau $  and a real constant $a\in [0,1]$, 
  we define {\em the Bouchaud trap model}, $\BTM(G,\boldsymbol \tau ,a)$ 
  as the Markov chain $X(t)$ on $\mathcal V$ whose dynamics is given by 
  \eqref{e:rates} and \eqref{e:generator}.
\end{definition}

In the original introduction of the model \cite{Bou92}, the trapping landscape 
$\boldsymbol \tau $ is given by a non-normalised Gibbs measure
\begin{equation}
  \boldsymbol \tau=\sum_{x\in \mathcal V} \tau_x\delta_x = 
  \sum_{x\in \mathcal V} e^{-\beta E_x}\delta_x,
\end{equation}
where $\beta >0$ is the inverse temperature and $E_x$ is seen as the 
energy at $x$. The rates $w_{xy}$ can be then expressed using the random 
variables $E_x$ instead of $\tau _x$,
\begin{equation}
    w_{xy}=  \nu \exp\big\{\beta \big((1-a)E_x-aE_y\big)\big\},
    \qquad\text{if $\<x,y\>\in \mathcal E$}.
\end{equation}

It is easy to check that $\boldsymbol \tau $ is a reversible measure 
for the Markov chain $X(t)$; 
the detailed balance condition is easily verified:
\begin{equation}
  \tau_x w_{xy}=\nu \tau_x^a \tau_y^a = \tau_y w_{yx}.
\end{equation}

Let us give some intuition about the dynamics of the chain $X$. 
Consider the embedded discrete-time Markov chain $Y(n)$: 
\begin{equation}
  X(t)=Y(n) \qquad \text{for all $S(n)\le t< S(n+1)$},
\end{equation}
 where $S(0)=0$ and $S(n)$ is the time of 
the $n^{\text{th}}$ jump of $X(\cdot)$. 

If $a=0$, the process $X$ is particularly simple. Its jumping rates 
$w_{xy}$ do not depend on the depth of the target vertex $y$.  $X$ 
waits at the vertex $x$ an exponentially distributed time with mean 
$\tau_x (\nu d_x)^{-1}$, where $d_x$ is the degree of $x$ in $G$. 
After this time, it jumps to one of the neighbours of $x$ chosen 
uniformly at random. Hence, the embedded discrete-time Markov chain 
$Y(n)$ is a simple random walk on the graph $G$ and $X(t)$ is its  
time change. The $a=0$ dynamics is sometimes referred to as \textit{Random 
  Hopping Times} (RHT) dynamics.

If $a\neq 0$, the jumping rates $w_{xy}$ depend on the target vertex. 
The process therefore does not jump uniformly to all neighbours of $x$. 
$Y(n)$ is no longer a simple random walk but a kind of  discrete 
Random Walk in Random Environment. To observe the effects of $a>0$ it 
is useful to consider a particular relatively deep trap $x$ with 
much shallower neighbours. In this case, as $a$ increases, the mean 
waiting time at $x$ decreases. On the other hand, if $X$ is located 
at some of the neighbours of $x$, then it is  attracted by the deep 
trap $x$ since $w_{yx}$ is relatively large. Hence, as $a$ increases, 
the process $X$ stays at $x$ a shorter time, but, after leaving it, it 
has larger probability to return there. We will see later that these 
two competing phenomena might exactly cancel in the long-time behaviour 
of well-chosen characteristics of the Markov chain.

We are not interested here in a natural line of questions which would 
be to find the best conditions under which the trapping mechanism is 
not crucial, and the BTM behaves as a simple random walk. On 
the contrary, we want to see how the trapping landscape can have a strong 
influence on the long time behaviour. Obviously, this can happen only 
if this trapping landscape is strongly inhomogeneous. 

Strong inhomogeneity can be easily achieved in the class of random 
landscapes with heavy tails, the essential hypothesis being that the 
expectation of the depth should be infinite. One of the assumptions 
we will use is therefore:
\begin{assumption}
  \label{a:stable}
  The depths $\tau_x$ are positive i.i.d.~random variables belonging to 
  the domain of the attraction of the totally asymmetric $\alpha $-stable 
  law with $\alpha \in (0,1)$. This means that there exists a slowly varying 
  function $L$ (i.e.,~for all $s>0$ $\lim_{u\to\infty} L(us)/L(u)=1$), 
  such that
  \begin{equation}
    \label{e:astable}
    \mathbb P[\tau_x\ge u]= u^{-\alpha }L(u).
  \end{equation}
\end{assumption}
Sometimes, to avoid unnecessary technical difficulties, we use the 
stronger assumption:
\begin{assumption}
  \label{a:limit}
  The depths $\tau_x$ are positive i.i.d.~random variables satisfying
  \begin{equation}
    \label{e:alimit}
    \lim_{u\to\infty} u^{\alpha }\mathbb P[\tau_x\ge u]= 1,
    \qquad \alpha \in (0,1).
  \end{equation}
\end{assumption}

These assumptions are satisfied at low temperature for the standard 
choice of the statistical physics literature:  for $-E_x$ being 
an i.i.d.~collection of exponentially distributed random variables with 
mean $1$. The depth of the traps then satisfies 
\begin{equation}
  \mathbb P[\tau_x\ge u]=\mathbb P[e^{-\beta E_x}\ge u]=u^{-1/\beta }, \qquad
  u\ge 1.
\end{equation}
Hence, if $\beta >1$, then Assumption \ref{a:limit} is satisfied  with 
$\alpha =1/\beta $, and the expected value of the depth diverges. 

\subsection{Examples of trap models}
\label{ss:examples}

Specifically, we will consider the following models:
\begin{enumerate}
  \item The Bouchaud model on $\mathbb Z$. Here $G$ is the integer 
  lattice, $\mathbb Z$, 
  with the nearest neighbour edges. The trapping landscape will be as in 
  Assumption~\ref{a:stable}. We will report in Section~\ref{s:onedim} 
  on work by \cite{FIN02, BC05, Cer06}.

  \item The BTM on $\mathbb Z^d$, $d>1$, with  a random landscape as 
  in Assumption~\ref{a:limit}. The main results about this case are 
  contained in \cite{BCM06,Cer03,BC06b}.

\end{enumerate}

We will also deal with a generalisation of the former setting, 
i.e.,~with a sequence of Bouchaud trap models, 
$\{\BTM(G_n,\boldsymbol \tau_n,a):n\in \mathbb N\}$. We will then consider 
different time scales depending on $n$ and write $X_n(t)$ for the 
Markov chain on the level $n$. In this case the law of $\tau_x$ may 
depend on $n$. Then the requirement $\mathbb E[\tau_x]=\infty$ is not 
necessary, as we will see.

\begin{enumerate}
  \addtocounter{enumi}{2}
  \item  The BTM in a large box in $\mathbb Z^d$ with periodic boundary 
  condition, here $G_n=\mathbb Z^d/n \mathbb Z^d$ is the torus of 
  size $n$ \cite{BC06}.

  \item The BTM on a large complete graph. Here we 
  consider a sequence of complete graphs with $n$ vertices. 
  This is the model that was originally proposed in \cite{Bou92}. 

  \item The Random Energy Model (REM) dynamics. We deal here with a 
  sequence of  $n$-dimensional hypercubes $G_n=\{-1,1\}^n$. The 
  landscape will be given by normally distributed (centred, with 
    variance $n$) energies $E_x$. This is the  case where $\tau$'s 
  are not heavy-tailed. 
\end{enumerate}

We will see that all the previous cases with exception of the 
one-dimensional lattice  behave very similarly. 

\subsection{Natural questions on trap models}

We will be interested in the long-time behaviour of the BTM. There 
are several natural questions to ask in order to quantify the 
influence of trapping. The most important question we will address in 
these notes is the question  of aging (at different time scales). Let 
us now define what we will call aging and subaging here. We will 
consider a time interval $[t_w,t_w+t]$, where the \textit{waiting 
  time} $t_w$ (or the age of the system) as well as the length $t$ of 
the time window (or duration of the observation) will grow to 
infinity. We will then consider various two-time functions, say 
$C(t_w,t_w+t)$, which depend on the trajectory of the Bouchaud trap 
model in the time window $[t_w,t_w+t]$. We will say that there is 
aging for these two-time functions iff
\begin{equation}
  \lim_{t_w\to \infty}
  C(t_w,t_w+\theta t_w)
\end{equation}
exists and is non-trivial. We will call this limit, $C(\theta )$, the 
aging function. We will say that there is subaging with exponent 
$\gamma <1$ iff
\begin{equation}
  \lim_{t_w\to \infty}
  C(t_w,t_w+\theta t_w^\gamma )
\end{equation}
exists and is again non-trivial.

We need now to define good two-time functions  in order to be able to 
deal with  aging for BTM. The following  functions are mostly studied:
\begin{itemize}
  \item[(a)] The probability that, conditionally on $\boldsymbol \tau $, 
  the process does not jump during the specified time interval $[t_w,t_w+t]$,
  \begin{equation}
    \Pi (t_w,t_w+t;\boldsymbol \tau )=\mathbb P\big[X(t')=X(t_w)\,
      \forall t'\in [t_w, t_w+t]\big|\boldsymbol \tau \big].
  \end{equation}

  \item[(b)] The probability that the system is in the same trap at both
  times $t_w$ and $t_w+t$,
  \begin{equation}
    R(t_w,t_w+t;\boldsymbol \tau )=\mathbb P[X(t_w)=X(t_w+t)|\boldsymbol \tau ].
  \end{equation}

  \item[(c)] And finally, the quantity
  \begin{equation}
    R^q(t_w,t_w+t;\boldsymbol \tau )=
    \mathbb E\Big[\sum_{x\in\mathcal V}
      [\mathbb P(X(t_w+t)=x|\boldsymbol \tau ,X(t_w))]^2\Big|\boldsymbol \tau \Big],
  \end{equation}
  which is the probability that two independent walkers will be at the 
  same site after time $t+t_w$ if they were at the same site at time 
  $t_w$, averaged over the distribution of the common starting point $X(t_w)$. 
\end{itemize}

These quantities are random objects: they still depend on the 
randomness of the trapping landscape $\boldsymbol \tau $. They are 
usually called {\em quenched} two-time functions. 
In addition to the quenched 
two-time functions, we will also consider their average over the random 
landscape. We define the {\em averaged} two-time functions:
\begin{equation}
  \begin{aligned}
    \Pi (t_w,t_w+t)&=\mathbb P\big[X(t')=X(t_w)\,
      \forall t'\in [t_w, t_w+t] \big],\\
    R(t_w,t+t_w )&=
    \mathbb P[X(t_w)=X(t_w+t) ],\\
    R^q(t_w,t_w+t )&=
    \mathbb E\Big[\sum_{x\in\mathcal V}
      [\mathbb P(X(t_w+t)=x|\boldsymbol \tau ,X(t_w))]^2 \Big].
  \end{aligned}
\end{equation}

We will state aging results for both quenched and averaged two-time 
functions for the various Bouchaud trap models given in 
Section~\ref{ss:examples}. We will strive to get the widest possible 
range of time scales where aging occurs. 

Even though our main motivation was the study of aging, there are 
many other questions about the long-time behaviour of the BTM which 
are of interest. For instance:
\begin{itemize}
  \item The behaviour of the environment seen from the particle. Here 
  the prominent feature of this environment seen from the position at 
  time $t$ is simply the depth of the trap $\tau_{X(t)}$ where the 
  process is located. We will give limit theorems for this quantity, 
  which are crucial for most of the aging results. 

  \item Nature of the spectrum of the Markov chain close to its edge. 
  Naturally, the long time behaviour of $X(t)$ can be understood from 
  the edge of the spectrum of the generator $L$. This question 
  deserves further study (see \cite{BF05,BF06} and also \cite{MB97}). 
  We intend to address this question for BTM in finite dimensions in 
  a forthcoming work.
  
  \item Anomalous diffusion. In the case where graph is $\mathbb Z^d$, 
  can we see that $X(t)$ is slow: for instance that 
  $\mathbb E[X(t)^2]\ll t$? Or get the tail behaviour of $|X(t)|$? 

  \item Scaling limit. Again in the case of $\mathbb Z^d$, is there a 
  scaling limit for the process $X(t)$, i.e.~a way to normalise space 
  and time so that 
  $X^\varepsilon (t)=\varepsilon X(t/h(\varepsilon ))$ converges to a 
  process on $\mathbb R^d$ which we can describe.
\end{itemize}

\subsection{References}
The physics literature on trap models is so abundant that we cannot 
try to be exhaustive. For earlier references on finite-dimensional 
questions see \cite{Mac85,BG90} where the anomalous character of the 
diffusion is given as well as a scaling limit. For aging questions on 
large complete graphs and relation to spin-glass dynamics see 
\cite{Bou92,BD95}, see also \cite{BCKM98} for a more global 
picture. For aging questions for the finite-dimensional model see 
\cite{MB96,RMB00,RMB01,BB03} among many other studies.

We will give references to mathematical papers at the end of every 
section, when necessary.


\section{The one-dimensional trap model}
\label{s:onedim}
We will consider in this section $\BTM(\mathbb Z,\boldsymbol \tau ,a)$ 
on the one-dimensional lattice, $\mathbb Z$, with nearest-neighbour 
edges. The random depths, $\tau_x$, will be taken to be i.i.d.,~in 
the domain of the attraction of an $\alpha $-stable law, $\alpha <1$, 
as in Assumption~\ref{a:stable}. There are several reasons why this 
particular graph should be treated apart. First, as usual, the 
one-dimensional model is easier to study. Second,   as we have 
already mentioned, the one-dimensional BTM  has some specific 
features that distinguish it from all other cases presented later in 
these notes.  

Another distinguishing feature (of technical character) is 
that, at present, the one-dimensional BTM is the only case where the 
asymmetric variant ($a>0$) has been rigorously studied. For technical 
convenience we choose here 
$\nu = \nu_a =  \mathbb E[\tau_0^{-a}]^2/2$, that is we set
\begin{equation}
  \label{ef:wij}
  w_{xy}=\frac 12\mathbb E[\tau_0^{-a}]^2 
  \tau_x^{-(1-a)}\tau_y^a,\qquad
    \text{if $|x-y|=1$.}
\end{equation}
We  set $X(0)=0$. We suppose that $\nu_a$ is finite for all 
$a\in[0,1]$, this is obviously the case if, e.g., $\tau_0>c$ a.s.~for 
some $c>0$. 

\subsection{The Fontes-Isopi-Newman singular diffusion}

The most useful feature of the one-dimensional BTM is that we can 
identify its scaling limit as an interesting one-dimensional singular 
diffusion in random environment introduced by Fontes, Isopi and 
Newman \cite{FIN02}. 

\begin{definition}[The F.I.N.~diffusion]
  \label{df:Z}
  Let $(x_i,v_i)$ be an inhomogeneous Poisson point process on 
  $\mathbb R\times(0,\infty)$ with intensity measure 
  $\d x\,\alpha v^{-1-\alpha }\d v$. Define the random discrete 
  measure $\rho = \sum_i v_i\delta_{x_i}$. We call $\rho $ the random 
  environment. Conditionally on~$\rho $, we define the {\em F.I.N. 
    diffusion}  $Z(s)$ as a diffusion process (with $Z(0)=0$) that 
  can be expressed as a time change of a standard one-dimensional 
  Brownian motion $B(t)$ with the speed measure $\rho$, as follows 
  \cite{IM65}: Denoting by $ \ell(t,y)$ the local time of the 
  standard Brownian motion $B(t)$ at $y$, we define
  \begin{equation}
    \phi_\rho (t)=\int_{\mathbb R} \ell(t,y)\rho (\d y)
  \end{equation}
  and its generalised right-continuous inverse
  \begin{equation}
    \psi_\rho (s)=\inf \{t>0:\phi_\rho (t)>s\}.
  \end{equation}
  Then $Z(s)=B(\psi_\rho (s))$. 
\end{definition}

The following proposition lists some of the properties of the 
diffusion $Z$ and the  measure $\rho $ that may be of interest. 

\begin{proposition}
  \label{p:Zproperties}
  \renewcommand\theenumi{(\roman{enumi})}
  \renewcommand\labelenumi{\theenumi}
  \begin{enumerate}
    \item The intensity measure of the Poisson point process 
    is non-integrable at $v=0$, therefore 
    the set of all atoms of $\rho $ is a.s.~dense in $\mathbb R$. 

    \item Conditionally on 
    $\rho $, the distribution of $Z(t)$, $t>0$,  is a 
    discrete probability measure $\nu^\rho_t=\sum_i w_i(t) \delta_{x_i}$, 
    with the same set of atoms as $\rho $.  

    \item The diffusion $Z$ has continuous sample paths.

    \item Define 
    $p_\rho (t,x_i)=\mathbb P[Z(t)=x_i|\rho ]/v_i$ for all atoms $x_i$ 
    of $\rho $ and $t>0$. The function $p_\rho $ 
    has a unique jointly continuous extension to $(0,\infty)\times \mathbb R$. 
    Moreover, $p_\rho $ satisfies the following  equation
    \begin{equation}
      \frac \partial {\partial t}p_\rho 
      (t,x)=\frac{\partial^2}{\partial \rho \partial x } p_\rho (t,x).
    \end{equation}
    The singular differential operator 
    $\partial^2 f /\partial\rho \partial x$ is defined (see, 
      e.g.,~\cite{DM76,KW82}) by 
    $h= \partial^2 f /\partial\rho \partial x$ if for some 
    $c_1,c_2\in \mathbb R$
    \begin{equation}
      f(x)=c_1+\int_0^x \Big(c_2 +\int_0^u h(v) \rho (\d v)\Big) \,\d u.
    \end{equation}

    \item  The diffusion $Z$ and its speed measure $\rho $ are 
    self-similar: for all $\lambda >0$, $t>0$ and $x\in \mathbb R$ 
    \begin{equation}
      \rho ([0,x])\laweq\lambda^{-1/\alpha }\rho ([0,\lambda x]) \quad 
      \text{and}\quad
      Z(t)\laweq\lambda^{-1}Z(t\lambda^{(1+\alpha )/\alpha }).
    \end{equation}
    Therefore, the diffusion $Z$ is anomalous.

    \item  There exist constants 
    $C$, $c$ such that for all $x$ and $t>0$
    \begin{equation}
      \mathbb P\big[|Z(t)|\ge x\big]\le C\exp\Big[
	-c\Big(\frac x{t^{\frac\alpha {1+\alpha}}}\Big)^{1+\alpha }\Big].
    \end{equation}
  \end{enumerate}
\end{proposition}
\begin{proof}
  Statement (i) is trivial, (ii) is proved in \cite{FIN02}. Claim 
  (iii) follows from (i), the continuity of sample paths of $B$, and 
  the definition of $Z$. (iv) is a non-trivial claim of the theory of 
  quasi-diffusions, see above references. The first part of (v) is a 
  direct consequence of the definition of $\rho $. The second part 
  then follows from the first one and from the well-known scaling 
  relations for the Brownian motion $B$ and its local time:
  \begin{equation}
    B(t)\laweq\lambda^{-1}B(\lambda^2 t)\qquad \text{and}\qquad
    \ell(t,x)\laweq\lambda^{-1}\ell(\lambda^2 t, \lambda x).
  \end{equation}
  The last claim is proved in \cite{Cer06}.
\end{proof}

\begin{remark}
  The scale of the upper bound in (vi) is probably optimal. The 
  corresponding lower bound was however never proved. The numerical 
  simulations and non-rigorous arguments  
  in \cite{BB03} however support this conjecture, and 
  give even exact values for constants $C$ and $c$. 
\end{remark}


\subsection{The scaling limit}
We now explain how the F.I.N.~diffusion appears as a scaling limit of 
the BTM. For all 
$\varepsilon \in (0,1)$ we consider the rescaled process 
\begin{equation}
  \label{e:defXe}
  X^\varepsilon (t)=\varepsilon X(t/\varepsilon c_\varepsilon ),
\end{equation}
where 
\begin{equation}
    \label{ef:ceps}
    c_\varepsilon =
    \big(\inf [t\ge0: \mathbb P(\tau_0>t)\le\varepsilon] \big)^{-1}.
\end{equation}
It follows from Assumption~\ref{a:stable} that there is a slowly 
varying function $L'$ such that 
$c_\varepsilon =\varepsilon^{1/\alpha }L'(1/\varepsilon )$. We will 
also consider the following rescaled landscapes,
\begin{equation}
  \label{e:taueps}
  \boldsymbol \tau^\varepsilon (\d x)=
  c_\varepsilon \sum_{y\in\mathbb Z}\tau_y\delta_{\varepsilon y}(\d x)=:
  \sum_{y\in \mathbb Z}\tau^\varepsilon_{\varepsilon y} 
  \delta_{\varepsilon y}(\d x).
\end{equation}

It is not so difficult to see that the distribution of 
$\boldsymbol \tau^\varepsilon $ converges to the distribution of 
$\rho $ as $\varepsilon \to 0$.
The next proposition states that it is possible to construct a 
coupling between different scales such that the convergence becomes 
almost sure.

\begin{proposition}[Existence of coupling]
  \label{p:coupling}
  There exists a family of measures 
  $\bar {\boldsymbol \tau}^\varepsilon $ and processes 
  $\bar X^\varepsilon $ constructed on the same probability space as  
  the measure $\rho $ and the Brownian motion $B$ such that 
  \renewcommand\theenumi{(\roman{enumi})}
  \renewcommand\labelenumi{\theenumi}
  \begin{enumerate}
    \item For all $\varepsilon >0$, 
    $\bar{ \boldsymbol \tau}^\varepsilon $  
    has the same distribution as $\boldsymbol \tau^\varepsilon $.

    \item The measures 
    $\bar {\boldsymbol \tau }^\varepsilon $ converge to $\rho $  
    vaguely and in the point-process sense, $\rho $-a.s.

    \item $\bar X^\varepsilon $ can be expressed as a time(-scale) 
    change of the Brownian motion $B$ with the speed measure 
    $\bar{\boldsymbol \tau}^\varepsilon $ (see Section~\ref{sss:tsc}). It has 
    the same distribution as 
    $X^\varepsilon $. 

  \end{enumerate}
\end{proposition}

We can now state the principal theorem of this section. It was proved 
in \cite{FIN02} for $a=0$ and in \cite{BC05} for $a>0$.

\begin{theorem}[Scaling limit of the one-dimensional BTM]
  \label{t:scalinglimit}
  As $\varepsilon \to 0$, for every fixed $t>0$ and all $a\in [0,1]$, 
  the distribution of 
  $\big(\bar X^\varepsilon (t), \bar \tau^\varepsilon_{\bar X^\varepsilon (t)}\big)$ 
  converges weakly and in the point-process sense to the distribution 
  of $\big(Z(t), \rho \big(\{Z(t)\}\big)\big)$, $\rho $-a.s.
\end{theorem}

The notion of convergence in the point-process sense used in the 
theorem was introduced in \cite{FIN02}. It is used here because we 
want to deal with quantities like 
$\mathbb P[X^\varepsilon(t)=X^\varepsilon (t')]$ to prove aging. The 
usual weak or vague convergences of measures are insensitive to such 
kind of quantities. This notion of convergence is defined by 

\begin{definition}[Point-process convergence]
  Given a family $\nu $, $\nu^\varepsilon $, $\varepsilon >0$, of 
  locally finite measures on $\mathbb R$,  
  we say that $\nu^\varepsilon $ converges {\em in
    the point process sense} to $\nu $, and write $\nu^\varepsilon
  \ppconv \nu $, as $\varepsilon \to 0$, provided the following 
  holds: If the atoms of $\nu $, $\nu^\varepsilon $ are,
  respectively, at the distinct locations $y_i$, $y_{i'}^\varepsilon$ 
  with weights $w_i$, $w_{i'}^\varepsilon$, then the subsets of
  $U^\varepsilon \equiv \cup_{i'}\{(y^\varepsilon_{i'},
      w^\varepsilon_{i'})\}$ of $\mathbb R\times(0,\infty)$ converge to
  $U\equiv\cup_i\{(y_i,w_i)\}$ as $\varepsilon \to 0$ in the sense
  that for any open $O$, whose closure is a compact subset
  of $\mathbb  R\times (0,\infty)$ such that its boundary contains no
  points of $U$, the number of points $|U^\varepsilon \cap O|$ in
  $U^\varepsilon \cap O$ is finite and equals $|U\cap O|$ for all
  $\varepsilon $ small enough.
\end{definition}

\begin{remark}
  The convergence of $\bar \tau^\varepsilon_{\bar X^\varepsilon (t)}$ 
  gives a description of the environment seen by the particle. More 
  explicitly, the distribution of the normalised depth of the trap 
  where $X$ is located at large time 
  $c_\varepsilon \tau_{X(t/\varepsilon c_\varepsilon) }$ converges 
  to the distribution of $\rho (\{Z(t)\})$.
\end{remark}

We now sketch the three main tools that are used in the proof of 
Theorem~\ref{t:scalinglimit}. 

\subsubsection{Time-scale change of Brownian motion} 
\label{sss:tsc}
To better 
understand how $(Z,\rho )$ arises as the scaling limit of 
$(X,\boldsymbol \tau )$, one should use the fact that not only 
diffusions, but also nearest-neighbour random  walks in dimension one,
can be expressed as time(-scale) change of the Brownian motion. The 
scale change is necessary only if $a\neq 0$, because the process 
$X(t)$ does not jump left or right with  equal probabilities.

We first define the time-scale change. Consider a locally-finite, 
discrete, non-random measure 
\begin{equation}
  \mu (\d x)=\sum_i w_i \delta_{y_i}(\d x),
\end{equation}
which has atoms with weights $w_i$ at positions $y_i$. The measure 
$\mu $ will be referred to as the {\em speed measure}. Let $S$ be a 
strictly increasing function defined on the set $\{y_i\}$. We call 
such $S$ the {\em scaling function}.  Let us introduce slightly nonstandard 
notation $\mm\mu S$ for the ``scaled measure'' 
\begin{equation}
  (\mm \mu S)(\d x)=\sum_i w_i\delta_{S(y_i)}(\d x).
\end{equation}
Similarly as in definition of $Z$, we define the function
\begin{equation}
  \label{ef:tch}
   \phi(\mu,S)(t)=\int_{\mathbb R} \ell(t,y)(\mm\mu S )(\d y)
\end{equation}
and the stopping time $\psi (\mu,S)(s)$ as the first time when 
$\phi (\mu,S) (t)=s$. The function $\phi (\mu,S)(t)$ is a nondecreasing, 
continuous function, and $\psi (\mu,S) (s)$ is its generalised right 
continuous inverse. It is an easy corollary of the results of \cite{Sto63} 
that the process
\begin{equation}
    \label{ef:defX}
  X(\mu ,S)(t):=S^{-1}(B(\psi (\mu,S) (t))) 
\end{equation}
is a  nearest-neighbour random walk on the set of atoms 
of $\mu $. Moreover, every nearest-neighbour random walk on a 
countable, nowhere-dense  subset of $\mathbb  R$ satisfying some mild 
conditions on transition probabilities can be expressed in this way.  
We call the process $ X(\mu ,S)$ the {\em time-scale change} of the 
Brownian motion. If $S=\Id$, the identity mapping, we speak only 
about the time change.

The following proposition summarises the properties of $X(\mu ,S)$ if 
the set of atoms of $\mu$ has no accumulation point. 
In this case we can suppose that the locations of atoms $y_i$ satisfy 
$y_i<y_j$ if $i<j$. 

\begin{proposition}[Stone, \cite{Sto63}]
    \label{pf:stone}
    The process $X(\mu ,S)(t)$ is a nearest-neighbour random walk on the 
    set $\{y_i\}$ of atoms of $\mu $. The waiting time in the state 
    $y_i$ is exponentially distributed with mean
    \begin{equation}
      2 w_i\frac{(S(y_{i+1})-S(y_i))(S(y_i)-S(y_{i-1}))}
      {S(y_{i+1})-S(y_{i-1})}.
    \end{equation}
    After leaving state $y_i$, $X(\mu ,S)$ enters states $y_{i-1}$ 
    and $y_{i+1}$ with respective probabilities 
    \begin{equation}
      \frac{S(y_{i+1})-S(y_i)}{S(y_{i+1})-S(y_{i-1})} \quad \text{and}\quad
      \frac{S(y_i)-S(y_{i-1})}{S(y_{i+1})-S(y_{i-1})}.
    \end{equation}
\end{proposition}

Using this proposition it is possible to express
the processes $X^\varepsilon (t)$ (see 
  \eqref{e:defXe}) as a time-scale change of the Brownian motion 
$B$. It is not surprise that 
$\boldsymbol \tau^\varepsilon  $  should be chosen as the speed 
measures. The scaling function is defined by
\begin{equation}
  \label{e:defS}
   S(x)=
  \begin{cases}
     \sum_{y=0}^{x-1}r_y,&\mbox{if $x \ge0,$}\\
     -\sum_{y=x}^{-1}r_y,&\mbox{otherwise}, 
  \end{cases}
\end{equation}
where
\begin{equation}
  \label{e:resist}
  r_x=\frac 12 \nu^{-1}_a \tau_x^{-a}\tau_{x+1}^{-a}.
\end{equation}
Observe that $\nu_a$ was chosen in such way that $\mathbb E[r_x]=1$.  
If $a=0$, $S$ is the identity mapping on $\mathbb Z$,  there is no 
scale change in this case.   Define further 
$S^\varepsilon(\cdot) = \varepsilon S(\varepsilon^{-1}\cdot)$. It is 
easy to check, using Proposition~\ref{pf:stone}, that the processes 
$X(\boldsymbol \tau^\varepsilon, S^\varepsilon)$ have the same 
distribution as $X^\varepsilon $. 

It is convenient to introduce  processes 
$\mathcal X^\varepsilon (t)$ that are the time change of the Brownian 
motion with speed measures  
$\mm{\boldsymbol \tau^\varepsilon}{ S^\varepsilon } $. Namely,
\begin{equation}
  \mathcal X^\varepsilon (t)=X(\mm {\boldsymbol \tau^\varepsilon } 
    {S^\varepsilon }, \Id)(t).
\end{equation}
The processes $\mathcal X^\varepsilon $  are related to 
$X^\varepsilon $ by
$ X^\varepsilon (t)=(S^\varepsilon)^{-1} (\mathcal X^\varepsilon (t))$. 

\subsubsection{Convergence of the fixed-time distributions} We have 
expressed the processes $\mathcal X^\varepsilon $ as the time change of 
the Brownian motion with the speed measure 
$\mm {\boldsymbol \tau^\varepsilon }{S^\varepsilon }$. We want 
to show that $\mathcal  X^\varepsilon $ and mainly $X^\varepsilon$ 
converge to $Z$. As stated in the following important theorem, it is 
sufficient to check the convergence of the speed measures to prove the 
convergence of fixed time distributions.  Observe that 
the theorem deals only with non-random measures.

\begin{theorem}[\cite{Sto63,FIN02}]
  \label{t:fin}
  Let $\mu^\varepsilon $, $\mu $ be a collection of {\em 
    non-ran\-dom} locally-finite measures, and let 
  $\mathcal Y^\varepsilon $, $\mathcal Y$ be defined by
  \begin{equation}
    \mathcal Y^\varepsilon (t)=X(\mu^\varepsilon  ,\Id)(t)\quad\text{and}\quad
    \mathcal Y(t)=X(\mu ,\Id)(t).
  \end{equation}
  For any deterministic $t_0>0$, let $\nu^\varepsilon$ denote the 
  distribution of $\mathcal Y^\varepsilon (t_0)$ and $\nu $ denote
  the distribution of $\mathcal Y(t_0)$. Suppose that
  \begin{equation}
    \mu^\varepsilon \vagconv \mu \quad \mbox{and} \quad 
    \mu^\varepsilon \ppconv  \mu \quad \mbox{as} \quad \varepsilon \to 0.
  \end{equation}
  Then, as $\varepsilon \to 0$,
  \begin{equation}
     \nu^\varepsilon \vagconv \nu \quad \mbox{and} \quad
     \nu^\varepsilon \ppconv  \nu.
  \end{equation}
  (Here $\vagconv$ stands for the vague convergence.)
\end{theorem}

\subsubsection{A coupling for walks on different scales.} The major 
pitfall of the preceding theorem is that it works only with sequences 
of deterministic speed measures. We want, however, to consider random 
speed measures $\boldsymbol \tau^\varepsilon$. As we have already 
remarked, it is not difficult to see that 
$\boldsymbol \tau^\varepsilon $ converge to  $\rho $ vaguely in 
distribution. However, it is not enough to make an application of 
Theorem~\ref{t:fin} possible. Here the coupling whose existence is 
stated in Proposition~\ref{p:coupling}(i) comes into play. It 
allows  to replace the convergence in distribution by the almost 
sure convergence. Then it is possible to apply Theorem~\ref{t:fin}.  
Let us construct this coupling.

Consider a two-sided L\'evy process (or $\alpha $-stable subordinator) $U(x)$, 
$x\in\mathbb R$, $U(0)=0$, with stationary and independent increments 
and cadlag paths defined by
\begin{equation}
  {\mathbb E}\big[e^{-\lambda (U(x+x_0)-U(x_0))}\big]=
   \exp\Big[x\alpha \int_0^\infty (e^{-\lambda w}-1)w^{-1-\alpha}\d w\Big].
\end{equation}
Let $\bar \rho $ be the random Lebesgue-Stieltjes measure on 
$\mathbb R$ associated to $U$, $\bar\rho (a,b]=U(b)-U(a)$. It is 
a known fact that $\bar\rho (\d x)=\sum_j v_j\delta_{x_j}(\d x)$, 
where $(x_j,v_j)$ is an inhomogeneous Poisson point process with 
density $\d x\,\alpha v^{-1-\alpha}\d v$, which means that $\bar\rho$ 
has the same distribution as $\rho$. 

For each fixed $\varepsilon >0$, we will now define the sequence of 
i.i.d.~random variables $\tau_x^\varepsilon $ such that 
$\tau_x^\varepsilon $'s are functions of $U$  and have the same 
distribution as $\tau_0$. 
Let  $G:[0,\infty)\mapsto[0,\infty)$ be such that 
\begin{equation}
  \label{ef:defG}
  {\mathbb P}(U(1)>G(x))=\mathbb P(\tau_0>x).
\end{equation}
It is well defined since $U(1)$ has continuous distribution, it is 
nondecreasing and right-continuous, and hence has nondecreasing 
right-continuous generalised inverse $G^{-1}$. 
\begin{lemma}
  \label{lf:cup}
   Let
   \begin{equation}
      \tau^\varepsilon_x:= 
      G^{-1} \big( 
	\varepsilon^{-1/\alpha}(U(\varepsilon(x+1))-U(\varepsilon x))\big).
   \end{equation}
   Then for any $\varepsilon >0$, the $\tau_x^\varepsilon $ are 
   i.i.d.~with the same law as $\tau_0$. 
\end{lemma}
\begin{proof}
  By stationarity and independence of increments of $U$ it is 
  sufficient to show 
  ${\mathbb P}(\tau_0^\varepsilon>t )=\mathbb P(\tau_0>t)$. However,
  \begin{equation}
    {\mathbb P}(\tau_0^\varepsilon>t )={\mathbb 
      P}(U(\varepsilon )>\varepsilon^{1/\alpha} G(t))
  \end{equation}
   by the definitions of $\tau^\varepsilon_0$ and $G$. The result 
   then follows from \eqref{ef:defG} and the scaling invariance of $U$: 
   $U(\varepsilon )\laweq \varepsilon^{1/\alpha} U(1)$. 
\end{proof}

Let us now define the random speed measures 
$\bar {\boldsymbol \tau} ^\varepsilon $ using the 
collections $\{\tau^\varepsilon_x\}$ from the previous lemma, 
\begin{equation}
  \label{ef:mubareps}
  \bar {\boldsymbol \tau }^\varepsilon (\d x)= 
  \sum_{i\in\mathbb Z}c_\varepsilon 
   \tau_i^\varepsilon \delta_{\varepsilon i}(\d x). 
\end{equation}
Finally, using $\tau^\varepsilon_x$ instead of $\tau_x$,
we define the scaling functions $\bar S^\varepsilon $ similarly as in 
\eqref{e:defS} and \eqref{e:resist}. The process $\bar X^\varepsilon $ 
is then given by 
$\bar X^\varepsilon = X(\bar{\boldsymbol \tau }^\varepsilon ,\bar S^\varepsilon )$, 
and the construction of the coupling from 
Proposition~\ref{p:coupling} is finished.

\subsubsection{Scaling limit}
\label{s:conv}

Using the three tools introduced above, we can now sketch the proof 
of Theorem~\ref{t:scalinglimit}. Actually, not many steps remain.

First, to prove the convergence of $\mathcal X^\varepsilon $ to $Z$ 
it is sufficient to verify the a.s.~convergence of the speed measures 
$\mm{\bar{\boldsymbol \tau }^\varepsilon}{\bar S^\varepsilon }$ to 
$\bar \rho $ and then apply Theorem~\ref{t:fin}. The proof of this 
convergence is not difficult, however, slightly lengthy. It can be 
found in \cite{FIN02} and \cite{BC05}.

\begin{proposition}
  \label{pf:conv}
  Let $\bar{\boldsymbol \tau }^\varepsilon$ and $\bar \rho$ be 
  defined as above. Then
  \begin{equation}
    \mm{\bar{\boldsymbol \tau }^\varepsilon}{\bar S^\varepsilon }
    \vagconv \bar \rho \quad \text{and}\quad
    \mm{\bar{\boldsymbol \tau }^\varepsilon}{\bar S^\varepsilon }
    \ppconv  \bar \rho \quad \text{as } 
    \varepsilon \to 0, \quad {\bar \rho }\mbox{-a.s.}
  \end{equation}
\end{proposition}

Finally, to pass from the convergence of $\mathcal X^\varepsilon $ to 
the convergence of $X^\varepsilon $ it is necessary to control the 
scaling functions $\bar S^\varepsilon$. 

\begin{lemma}
  \label{lf:unif}
   As $\varepsilon \to 0$ we have 
   \begin{equation}
     \bar S^\varepsilon(\varepsilon \lfloor \varepsilon^{-1} y\rfloor )\to 
     y, \qquad  {\bar \rho }\mbox{-a.s.},
   \end{equation}
   uniformly on compact intervals.
\end{lemma}

Observe that this lemma also implies that the embedded discrete-time 
random walk $Y$ converges, after a renormalisation, to the Brownian 
motion, independently of the value of $a$. This is valid also if $a>0$ and 
the discrete-time embedded process $Y$ is not a simple random walk 
but a random walk in random environment. 

Since Lemma~\ref{lf:unif} is one of the key 
parts of the proof of the scaling limit for $a\neq 0$ we prove it here. 

\begin{proof}[Proof of Lemma~\ref{lf:unif}]
  We consider only $y>0$. The proof for $y<0$ is very similar. By 
  definition of $\bar S^\varepsilon $ we have 
  $\varepsilon \bar S^\varepsilon (\lfloor \varepsilon^{-1} y\rfloor )= 
  \varepsilon \sum_{j=0}^{\lfloor \varepsilon^{-1} 
    y\rfloor-1}\bar r_j^\varepsilon $, 
  where for fixed $\varepsilon$ the sequence $\bar r_i^\varepsilon $ is an 
  ergodic sequence of bounded positive random variables. Moreover, 
  $\bar r_i^\varepsilon $ is independent of all $\bar r_j^\varepsilon $ with 
  $j\notin\{i-1,i,i+1\}$. The $\bar\rho $-a.s.\  convergence for 
  fixed $y$ is then a consequence of the strong law of large numbers 
  for triangular arrays.  Note that this law of large numbers can be 
  easily proved in our context using the standard 
  methods, because the variables $\bar r_i^\varepsilon $ are bounded and 
  thus their moments of arbitrary large degree are finite.  The uniform 
  convergence on compact intervals is easy to prove using the fact
  that $\bar S^\varepsilon $ is increasing and the identity function is continuous. 
\end{proof}

\subsection{Aging results}
The aging results for the one-dimensional BTM follow essentially 
from Theorem~\ref{t:scalinglimit}. To control the two-time functions $R$ 
and $R^q$ it is only necessary to extend its validity to the joint 
distribution of $(X^\varepsilon (1), X^\varepsilon (1+\theta ))$ at 
two fixed times, which is not difficult. This extension then yields 
the following aging result. 

\begin{theorem}[Aging in the one-dimensional BTM]
  \label{TF:AGING}
  For any $\alpha \in (0,1)$, $\theta >0$ and $\inp\in[0,1]$ there exist 
  aging 
  functions $ R_1(\theta)$, $R^q(\theta)$ such that
  \begin{equation}
    \begin{split}
      \lim_{t_w\to\infty} R (t_w,t_w+\theta t_w) &=
      \lim_{t_w\to\infty} \mathbb  E \mathbb P[X((1+\theta)t_w) = 
	X(t_w)|\boldsymbol \tau ]=R_1 (\theta),\\
      \lim_{t_w\to\infty} R^q(t_w,t_w+\theta t_w) &=
      \lim_{t_w\to\infty} \mathbb  E 
      \sum_{i\in \mathbb Z}
      [\mathbb P(X((1+\theta)t_w)=i|\boldsymbol \tau ,X(t_w))]^2=R^q(\theta).
    \end{split}
  \end{equation}
  Moreover, $R_1 (\theta)$ and $R^q(\theta)$ can be expressed using the 
  analogous quantities defined using the singular diffusion $Z$:
  \begin{equation}
    \begin{aligned}
      R_1 (\theta)&=\mathbb E\mathbb P[Z(1+\theta)=Z(1)|\rho],\\
      R^q(\theta)&=\mathbb E\sum_{x\in \mathbb R}
      [\mathbb P(Z(1+\theta)=x|\rho ,Z(1))]^2.
     \end{aligned}
  \end{equation}
\end{theorem}

\begin{remark}
  1.  This result is contained in \cite{FIN02} for $\inp=0$ and 
  in \cite{BC05} for $a>0$.  

  2. Let us emphasise that the functions $R_1(\theta )$, 
  $R^q(\theta )$ do not depend on the parameter $a$, since the 
  diffusion $Z(t)$ and the measure $\rho $ do not depend on it.  This 
  is the result of the compensation of shorter visits of deep traps 
  by the attraction to them.

  3. It should be also underlined that only averaged functions are 
  considered in the theorem. For a fixed realisation of 
  $\boldsymbol \tau $  there is no limit of 
  $R(t_w,t_w+\theta t_w;\boldsymbol \tau )$. From the proof of 
  Theorem~\ref{TF:AGING}, it is however not difficult to derive the 
  following weaker result (see Theorem~1.3 in~\cite{Cer06}).
\end{remark}

\begin{theorem}[Quenched aging on $\mathbb Z$, in distribution]
  \label{t:onedquen}
  As $t_w\to\infty$, the distribution of 
  $R(t_w,t_w+\theta t_w;\boldsymbol \tau )$ converges weakly to the 
  distribution of $\mathbb P[Z(1+\theta )=Z(1)|\rho ]$.
\end{theorem}

\subsection{Subaging results}
In the case of the two-time function $\Pi $ much shorter times $t$ 
should be considered, $t\ll t_w$, i.e.~subaging takes place: 

\begin{theorem}
  \label{tf:subaging}
  For any $\alpha \in (0,1)$, $\theta >0$ and $\inp\in[0,1]$ there exist an
  aging function $\Pi_{1,a} (\theta)$ such that
  \begin{multline}
      \lim_{t\to\infty}\Pi (t,t+f_a(t,\theta ))=
      \lim_{t\to\infty}\mathbb E\mathbb 
      P\big[X(t')=X(t)\forall t'\in[t,t+f_a(t,\theta )]|\boldsymbol 
	\tau \big]=
      \Pi_{1,a} (\theta),
  \end{multline}
  where the function $f_a$ is given by 
  \begin{equation}
    f_a(t,\theta)=\theta t^{\gamma (1-\inp)} L(t)^{1-\inp}.
  \end{equation}
  Here we use $\gamma $ to denote  the {\em subaging exponent}, 
  $\gamma = (1+\alpha )^{-1}$,
  and $L(t)$ is a slowly varying function that is determined only by the 
  distribution of $\tau_0$.   
  The function $\Pi_{1,a} (\theta)$ can again be written using the singular 
  diffusion $Z$,
  \begin{equation}
    \label{ef:thm2}
    \Pi_{1,a} (\theta)=\int_0^\infty g_a^2(\theta u^{a-1})\d F(u),
  \end{equation}
  where $F(u)=\mathbb E \mathbb P[\rho(\{Z(1)\})\le u| \rho ]$, and where 
  $g_a(\lambda )$ is the Laplace transform  
  of the random variable $\nu_a \tau_0^a$, 
  \begin{equation}
    g_a(\lambda )=\mathbb E(e^{-\lambda \nu_a \tau_0^a}).  
  \end{equation}
  If $\inp=0$, \eqref{ef:thm2} can be written as
  \begin{equation}
    \Pi_{1,0} (\theta)=\int_0^\infty e^{-\theta/u}\d F(u).
  \end{equation}
\end{theorem}

\begin{remark}
  1. As can be seen, in this case the function $\Pi_{1,a} (\theta)$ 
  depends on $\inp$. This is not surprising since the compensation by 
  attraction has no influence here and the jump rates clearly depend 
  on $\inp$.  

  2. Of course, an analogous result to Theorem~\ref{t:onedquen} holds 
  for quenched subaging in distribution.
\end{remark}

In the RHT case, $a=0$, the proof of Theorem~\ref{tf:subaging} is 
straightforward. It follows from the convergence of the distribution 
of $\tau^\varepsilon_{X^\varepsilon (1)}$ to the distribution of 
$\rho(\{Z(1)\})$ as stated in Theorem~\ref{t:scalinglimit}. This then 
implies the convergence of the distribution of 
$\tau_{X(t_w)}/(t_w^{\gamma }L(t_w))$. 

The proof in the case $a>0$ is more complicated. It essentially 
involves a control of the distribution of the depth of those traps 
that are nearest neighbours of the traps where $X(t_w)$ is with a 
large probability. It turns out that this distribution converges to 
the the distribution of $\tau_0$, i.e.~there is nothing special on 
the neighbours of deep traps.

\medskip

The behaviour  of the two-point functions $\Pi (t_w,t+t_w)$ and 
$R (t_w,t+t_w)$ is not difficult to understand and guess. We give 
here a heuristic explanation for these results, first in the case 
$a=0$. After the first $n$ jumps the process typically visits 
$O(n^{1/2})$ sites. The deepest trap that it finds during $n$ jumps 
has therefore a depth of order $O(n^{1/2\alpha })$, which is the 
order of the maximum of $n^{1/2}$ heavy-tailed random variables 
$\tau_x$ as can be verified from Assumption~\ref{a:limit}. This trap 
is typically visited $O( n^{1/2})$ times. Since the depths are in the 
domain of attraction of an $\alpha $-stable law with $\alpha <1$, the 
time needed for $n$ jumps is essentially determined by the time spent 
in the deepest trap. This time is therefore 
$O(n^{(1+\alpha )/2\alpha })$. Inverting this expression we get that 
the process visits typically $O(t^{\alpha \gamma })$ sites before 
time $t$. The deepest traps it finds during this time have a depth of 
order $t^\gamma $. Moreover, the process  is located in one of these 
deep traps at time $t$. From this, one sees that  the main 
contribution to  quantity $R(t_w,t_w+t)$  comes from the trajectories 
of $X$ that, between times $t_w$ and $t_w+t$, leave  the original 
site $X(t_w)$ a number of times of order $t_w^{\alpha\gamma  }$, and 
then return to it. Each visit of the original site lasts an amount of 
time of order $t_w^{\gamma } $, which is the time scale on which 
$\Pi $ ages. 

For $a>0$, such heuristics is not directly accessible. However, the 
Theorem~\ref{t:scalinglimit} yields that $\tau_{X(t_w)}$ is of the 
same order, $O(t_w^\gamma )$, as in the RHT case. Each visit of 
$X(t_w)$ lasts a shorter time, $O(t_w^{\gamma (1-a)})$, as follows 
from the definition of the process and the fact that the depths of 
the neighbours of $X(t_w)$ are $O(1)$. On the other hand, the process 
makes more excursions from  $X(t_w)$, their number being 
$O(t_w^{\gamma (\alpha +a)})$. 

\subsection{Behaviour of the aging functions on different time scales}
\label{ss:differen}

Having found two interesting time scales $t=O(t_w)$ and 
$t=O(t_w^\gamma )$ in the model%
\footnote{We suppose here that $a=0$ and Assumption~\ref{a:limit} holds, that 
  is $L(t)\to 1$.}
one may ask if there are other interesting time scales 
for the functions $\Pi $ and $R$. This question was raised by Bouchaud and 
Bertin in \cite{BB03}. The negative answer was given in \cite{Cer06}: 

\begin{theorem}[Behaviour of $\Pi $ on different scales]
  \label{t:shortscales}
  \renewcommand\theenumi{(\alph{enumi})}
  \renewcommand\labelenumi{\theenumi}
  \begin{enumerate}
    \item {\em Short time scales.~~} Let $f(t)$ be an 
    increasing function satisfying $t^\kappa \ge f(t)\ge t^\mu $ 
    for all $t$ large and for some $\gamma > \kappa \ge \mu >0$. Then
    \begin{equation}
      \label{e:Pist}
      \lim_{t\to\infty}
      \Big(\frac{f(t)}{t^\gamma }\Big)^{\alpha -1} \big(1-\Pi (t,t+f(t) )\big)=K_1,
    \end{equation}
    with $0 < K_1 <\infty$.

    \item 
    {\em Long time scales.} Let $g(t)$ be such that $t^\gamma =o(g(t))$. 
    Then
    \begin{equation}
      \label{e:Pilt}
      \lim_{t\to\infty}
      \Big(\frac{g(t)}{t^\gamma }\Big)^{\alpha } \Pi \big(t,t+g(t) \big)=K_2,
    \end{equation}
    with $0 < K_2 <\infty$.

    \item {\em Behaviour of $\Pi_{1,0} (\theta )$.} 
    The function $\Pi (\theta )$ satisfies
      \begin{gather}
	\label{e:Pias}
	\lim_{\theta \to 0} \theta^{\alpha -1}(1-\Pi_{1,0} (\theta ))=K_1,\\
	\label{e:Piaslt}
	\lim_{\theta \to \infty} \theta^{\alpha }\Pi_{1,0} (\theta )=K_2.
      \end{gather}
  \end{enumerate}
\end{theorem}

\begin{remark}
  We emphasise that the constant $K_1$ occur both in \eqref{e:Pist} 
  and \eqref{e:Pias}. That means that the behaviour of 
  $\Pi_{1,0}(\theta )$ at $\theta \sim 0$ gives also the behaviour of 
  $\Pi (t, t+f(t))$ for $f(t)\ll t^\gamma $. An analogous remark 
  applies for time scales $f(t)\gg t^{\gamma }$. Both constants $K_1$ and 
  $K_2$ can again be expressed using the F.I.N.~diffusion. The 
  formulae, and also similar (slightly weaker) results for the 
  two-time function $R$ can be found in \cite{Cer06}. 
\end{remark}

We give again a heuristic description of this result. As we know, at 
time $t$ the process $X$ is typically in a trap  of depth 
$O(t^{\gamma })$, it needs a time of the same order to jump out. In 
Theorem \ref{t:shortscales}(a) we look at $1-\Pi (t,t+f(t) )$ with 
$f(t) \ll t^\gamma $, that is at the probability that a jump occurs 
in a time much shorter than $t^{\gamma }$.  There are essentially 
two possible extreme strategies which lead to such an event:
\begin{enumerate}
  \item $\tau_{X(t)}$ has the typical order $t^\gamma $ but the 
  jump occurs in an exceptionally short time.

  \item $X(t)$ is in an atypically shallow trap
  and stays there a typical time.
\end{enumerate}
In \cite{Cer06} it is proved that the second strategy dominates. 
Therefore, one has to study the probability of being in a very 
shallow trap or, equivalently, to describe the tail of 
$\mathbb P[\tau_{X(t)}/t^\gamma  \le u]$ for $u\sim 0$. To control 
this tail  the proof use the fact that although the BTM never 
reaches equilibrium in a finite time in infinite volume, it is nearby  
equilibrium if one observes only traps that are much shallower than 
the typical depth $t^\gamma $ on intervals that are small with 
respect to the typical size of $X(t)$.  This puts on a rigorous 
basis, at least in dimension one,  the concept of {\em local 
  equilibrium} that was introduced in the physics literature by 
\cite{RMB00}. The concept does not give the right predictions for the 
values of the limiting functions $R_1(\theta )$ and 
$\Pi_{1,0} (\theta )$ but it is useful to describe their asymptotic 
behaviour. 
A very similar heuristics applies also for time scales $f(t)\gg t^\gamma $.

\subsection{References} In the RHT case, $a=0$, this model has first 
time been studied  by Fontes, Isopi and Newman in \cite{FIN99}.  It 
was used there as a tool to control the behaviour of the voter model 
with random rates. It was proved there that the process is 
sub-diffusive, i.e.~$\mathbb E[X(t)/\sqrt t]\to 0$ as $t\to\infty$, 
and that the dynamics localises in the sense that (for a.e.~$\boldsymbol \tau $)
\begin{equation}
  \sup_{x\in \mathbb Z}\mathbb P[X(t)=x|\boldsymbol \tau ]\not\to 0 \qquad 
  \text{as $t\to\infty$}.
\end{equation}
That means that there will be always a site (dependent on time $t$) where 
the process $X$ can be found with  a non-negligible probability. This 
localisation occurs only in one-dimensional BTM and is at the heart of 
the majority of the differences between $\mathbb Z$ and other graphs. 

The scaling limit has been established for the case $a=0$ in 
\cite{FIN02} and by \cite{Cer03,BC05} for $a\neq 0$. Following 
\cite{BB03}, the results of Section~\ref{ss:differen} have been given 
by \cite{Cer06}. 

For aging in another interesting one-dimensional  dynamics, 
i.e.~Sinai's Random Walk, see \cite{DGZ01,BF06}. 
Let us mention two interesting open questions on this 
one-dimensional trap model:

(a) What is the behaviour of the edge of the spectrum for the generator 
of the dynamics. This might be close to, but easier than the same 
question solved for  Sinai's Random Walk by \cite{BF06}. 

(b) What is the influence of a drift in the BTM? Monthus \cite{Mon04} 
gives a very interesting picture based on renormalisation arguments.


\section{The trap model in dimension larger than one}
\label{s:ddim}
After resolving the BTM  on $\mathbb Z$, the next natural step is to 
study the Bouchaud model on the  $d$-dimensional lattice, 
$\BTM(\mathbb Z^d,\boldsymbol \tau ,0)$, $d>1$. Observe that we set 
$a=0$, that means that only the RHT dynamics is considered. In this 
section we always assume that Assumption~\ref{a:limit} holds.

\subsection{The fractional-kinetics process}

As for the one-dimensional model, we first  identify 
a scaling limit of the BTM on $\mathbb Z^d$. The result of this section 
is contained in the 
forthcoming paper~\cite{BC06b}. We will, from now on, use frequently the 
theory of L\'evy processes and subordinators. A very short summary of 
this theory can be found in Appendix~\ref{ap:subord}.

Let us first define the process that appears as the scaling limit.

\begin{definition}[Fractional kinetics]
  Let $B_d(t)$ be the standard $d$-dimensional Brownian motion started at 
  $0$ and let $V$ be the $\alpha $-stable subordinator given by its 
  Laplace transform 
  $\mathbb E[e^{-\lambda V(t)}]=e^{-t\lambda^\alpha }$. Let 
  $T(s)=\inf\{t:V(t)>s\}$  be the inverse of $V(t)$. 
  We define the \textit{fractional-kinetics process} $\Psi_d$ by  
  \begin{equation}
    \Psi_d (s)=B_d(T(s)).
  \end{equation}
\end{definition}  

We list here without proofs several properties of the process $\Psi_d $.

\begin{proposition}
  \label{p:fk}
  1. The stable subordinator $V$ is strictly increasing, therefore its 
  inverse $T$ and also the process $\Psi_d$ are continuous. 

  2. The name of the process is due to the following fact. Let 
  $p(t,\cdot)$ be the probability density of $\Psi_d (t)$. Then $p$ is a 
  solution of the fractional kinetic equation,
  \begin{equation}
    \label{e:fract}
    \frac{\partial^\alpha }{\partial t^\alpha } p(t,x)=
    \frac 12 \Delta p(t,x) + \delta (0)\frac {t^{-\alpha }}{\Gamma 
      (1-\alpha )}.
  \end{equation}
  Here, the fractional derivative $\partial^\alpha p(t,x) /\partial t^\alpha $ is 
  the inverse Laplace transform of $s^\alpha \tilde p(s,x)$, where 
  $\tilde p(s,x)=\int_{0}^\infty e^{-st} p(t,x) \,\d t$ is the usual 
  Laplace transform. The equation \eqref{e:fract} should be understood 
  in the weak sense, i.e.~it holds after the integration against smooth 
  test functions.

  3. The process $\Psi_d $ is not Markov, as can be seen easily from the 
  previous point. 

  4. The fixed-time distribution of $\Psi_d$ is the Mittag-Leffler distribution,
  \begin{equation}
    \mathbb E\big(e^{i \xi \cdot \Psi_d(t)}\big)=E_\alpha (-|\xi |^2 t^\alpha ),
  \end{equation}
  where $E_\alpha (z)=\sum_{m=0}^\infty z^m/\Gamma (1+m\alpha )$.

  5. The process $\Psi_d$ is self-similar:
  \begin{equation}
    \Psi_d(t)\laweq \lambda^{-\alpha /2}\Psi (\lambda t).
  \end{equation}
\end{proposition}

The process $\Psi_d$ is well known in the physics literature (see 
  \cite{Zas02} for a broad survey and earlier references). It is 
the scaling limit of a very classical object, a Continuous Time 
Random Walk (CTRW) introduced by \cite{MW65}. More precisely 
consider a simple random walk $Y$ on $\mathbb Z^d$ and a sequence of 
positive i.i.d.~random variables $\{s_i:i \in \mathbb N\}$ with 
the distribution in the domain of attraction of an $\alpha $-stable law. 
Define the CTRW $U(t)$ by 
\begin{equation}
  U(t)=Y(k) \qquad \text{if } t \in \Big[\sum_{i=1}^{k-1} s_i, 
    \sum_{i=1}^{k} s_i\Big).
\end{equation}
It is proved in \cite{MS04} that there is a constant 
$C$ (depending only on the distribution of $s_i$) such that
\begin{equation}
  C n^{-\alpha /2} U(tn)\xrightarrow{n\to\infty}\Psi_d(t).
\end{equation}

\subsection{Scaling limit}

Observe that  
the Bouchaud model (for $a=0$) can be expressed as a time change 
of the simple random walk. The time-change process is crucial for us:
\begin{definition}
  \label{d:clock}
  Let $S(0)=0$ and let $S(k)$, $k\in \mathbb N$, be the time of the  
  $k^{\text{th}}$ jump of $X$. For $s \in \mathbb R$ we define 
  $S(s)=S(\lfloor s \rfloor)$. We call $S(s)$ the \textit{clock process}.
  Obviously, $X(t)=Y(k)$  for all $S(k)\le t< S(k+1)$.
\end{definition}

The following result shows that the limit of the $d$-dimensional 
Bouchaud model and its clock process on $\mathbb Z^d$ ($d \ge 2$) is 
trivial, in the sense that it is identical with the scaling limit of 
the much simpler (``completely annealed'') dynamics of the CTRW. 

\begin{theorem}[Scaling limit of BTM on $\mathbb Z^d$]
  \label{t:dscaling}
  Let
  \begin{equation}
    f(n)=\begin{cases}
      n^{\alpha /2}(\log n)^{(1-\alpha )/2},&\text{if $d=2$,}\\
      n^{\alpha /2},&\text{if $d\ge 3$}.
    \end{cases}
  \end{equation}
  Then for all $d\ge 2$ and for a.e.~$\boldsymbol \tau $,
  \begin{equation}
    \frac {C_d(\alpha ) X(nt)}{f(n)}
    \xrightarrow{n\to\infty} \Psi_d (t) \quad \text{and} \quad
    \frac {S\big( C_d(\alpha )^{-2}f(n)^2 s\big)}n
    \xrightarrow{n \to\infty} V(s).
  \end{equation}
  weakly in the Skorokhod topology on $D([0,T],\mathbb R^d)$ (the space 
    of cadlag functions from $[0,T]$ to $\mathbb R^d$). If  
  $G_d(0)$ denotes  Green's function of the $d$-dimensional random walk 
  at~$0$, then
  \begin{equation}
    C_d(\alpha )=
    \begin{cases}
      \sqrt{  2 \pi^{1-\alpha } \Gamma (1-\alpha )\Gamma (1+\alpha )} 
      ,&\text{if $d=2$,}\\
      \sqrt{ d G_d(0)^\alpha \Gamma (1-\alpha )\Gamma (1+\alpha )}
      ,&\text{if $d\ge 3$}.
    \end{cases}
  \end{equation}
\end{theorem}

The main ideas of the proof of this theorem will be explained in 
Section~\ref{ss:coarse}. At this place, let us only compare the 
fractional-kinetics process $\Psi_d$ with the F.I.N.~diffusion~$Z$. 
Both these processes are defined as a time change of the Brownian 
motion $B_d(t)$. The clock processes however differ considerably. For 
$d=1$, the clock equals $\phi (t)=\int \ell(t,x)\rho (\d x)$, where 
$\rho $ is the random speed measure obtained as the scaling limit of 
the environment. Moreover, since $\ell$ is the local time of the 
Brownian motion $B_1$, the processes $B_1$ and $\phi $ are {\em 
  dependent}. For $d\ge 2$, the Brownian motion $B_d$ and the clock 
process, i.e.~the stable subordinator $V$, are \textit{independent}. 
The asymptotic independence of the clock process $S$ and the location 
$Y$ of the BTM is a very remarkable feature distinguishing $d\ge 2$ 
and $d=1$. It explains the ``triviality'' of the scaling limit in 
dimension $d\ge 2$, but  is, by no means, trivial matter to prove. We 
will come back to an intuitive explanation of the independence in 
Section~\ref{ss:verif}.  Note also that nothing like a scaling limit 
of the random environment appears in the definition of $\Psi_d$, 
moreover, the convergence holds $\boldsymbol \tau $-a.s. The absence 
of the scaling limit of the environment in the definition of $\Psi_d $ 
transforms into the non-Markovianity of $\Psi_d $. Note however that 
it is considerably easier to control the behaviour of $\Psi_d $ than 
of $Z$ even if $\Psi_d $ is not Markov: many quantities related to 
$\Psi_d$ can be computed explicitly, as can be seen from 
Proposition~\ref{p:fk}.

\subsection{Aging results}

The following two theorems describe the aging behaviour of the 
two-time functions $R$ and~$\Pi $. 

\begin{theorem}[Quenched aging on $\mathbb Z^d$]
  \label{t:aging}
  For all $\alpha \in (0,1)$ and  $d\ge 2$  there exists a 
  deterministic function $R(\theta )$  independent of $d$ (but 
    dependent on $\alpha $)
  such that for $\mathbb P$-a.e.\ 
  realisation of the random environment $\boldsymbol \tau  $
  \begin{equation}
    \lim_{t_w\to\infty}R (t_w,t_w+\theta t_w;\boldsymbol \tau )=R (\theta ).
  \end{equation}
  The function $R(\cdot)$ can be written explicitly:  
  Let $\Asl_\alpha (u)$ be the distribution function of the 
  generalised arcsine law with parameter $\alpha $,
  \begin{equation}
    \Asl_\alpha (u):=
    \frac{\sin\alpha \pi }{\pi }
    \int_0^u u^{\alpha -1}(1-u)^{-\alpha }\,\d u.
  \end{equation}
  Then $R(\theta )=\Asl_\alpha (1/1+\theta )$.
\end{theorem}

\begin{theorem}[Quenched (sub-)aging on $\mathbb Z^d$]
  \label{t:subaging}
  For all $\alpha \in (0,1)$ and  $d\ge2$ there exists a 
    deterministic function $\Pi_d (\theta )$ such that for 
  $\mathbb P$-a.e.\ realisation of the random environment~$\boldsymbol \tau  $
  \begin{equation}
    \lim_{t_w\to\infty}
    \Pi \big(t_w,t_w+\theta f(t_w);\boldsymbol \tau \big)=\Pi_d (\theta ),
  \end{equation}
  with 
  \begin{equation}
    f(t_w)=
    \begin{cases}
      \frac {t_w}{\log t_w},&\text{if $d=2$,}\\
      t_w,&\text{if $d\ge 3$.}
    \end{cases}
  \end{equation}
  The function $\Pi_d$ does depend on $d$. Moreover, for all 
  $\theta >0$ it satisfies 
  \begin{equation}
    \lim_{d\to\infty}\Pi_d(\theta )=R(\theta ).
  \end{equation}
\end{theorem}

\begin{remark} 
  1. Both theorems are proved in \cite{BCM06} for $d=2$ and in 
  \cite{Cer03} for $d\ge 3$. The proofs are  relatively technical and 
  exceed the scope of these notes. The main ideas however do not use 
  specific properties of the integer lattice $\mathbb Z^d$ and can be 
  generalised to different graphs. These ideas will be explained in 
  Section~\ref{s:general}. The other important part of the proof, 
  that is the coarse-graining of the trajectory of the process, is 
  explained in the next subsection.  

  2. The function  $\Pi_d$ can  also be explicitly calculated but the 
  formula is tedious \cite{BCM06,Cer03}.

  3.  We will see later that the function $R$ is closely related to 
  the arcsine law for L\'evy processes. As we will also see,  the same 
  function appears as the limit  in the case of the BTM on a large 
  complete graph. Therefore, the mean-field dimension of the BTM for 
  the two-time function $R$ is $d=2$. 

  4. Both presented results are {\em quenched}, i.e.~they hold for 
  a.e.~$\boldsymbol \tau $. This should be compared with the results 
  for the $d=1$ case (Theorems~\ref{TF:AGING}, \ref{t:onedquen} and 
    \ref{tf:subaging}) where only averaged aging holds. The reason 
  for this difference is explained in Section~\ref{ss:verif} below. 
  It is, of course, trivial to get averaged results from 
  Theorems~\ref{t:aging} and \ref{t:subaging}: the dominated 
  convergence theorem yields
  \begin{equation}
    \begin{aligned}
     & \lim_{t_w\to\infty}R (t_w,t_w+\theta t_w)=R (\theta ),\\
     & \lim_{t_w\to\infty}
      \Pi \big(t_w,t_w+\theta f(t_w)\big)=\Pi (\theta ).
  \end{aligned}
  \end{equation}

  5. The fact that in $d\ge 3$ the time scale $f(t_w)=t_w$ is the same 
  for both functions $R$ and $\Pi $ is a consequence of the 
  transience of the simple random walk. The traps are  visited only a 
  finite number of times, so that different time scales cannot appear 
  for $R$ and~$\Pi $.
\end{remark}

\subsection{The coarse-graining procedure}
\label{ss:coarse}

We would like to describe here the coarse-graining procedure which 
was introduced in \cite{BCM06}, and which is the main tool in proving 
Theorems~\ref{t:dscaling}--\ref{t:subaging}.  Even this short sketch 
of the procedure might be considered technical and can be skipped on 
a first reading. The reader can also decide to return here after 
being acquainted  with the general ideas of 
Sections~\ref{ss:universal} and~\ref{ss:verif}. We will deal here 
only with the convergence of the clock process to an $\alpha $-stable 
subordinator. For the sake of concreteness we set $d=2$, however 
the same arguments apply also for $d>2$. The discussion in 
this section is valid \textit{for a.e.~realisation of the random 
  environment $\boldsymbol \tau $}.  

The random environment is heavy tailed. Therefore, as for $d=1$, the 
behaviour of the clock process $S(k)$ is determined by the time spent in 
the deepest traps that the process visits during the first $k$ steps. It 
is thus necessary to find the depth scale of these traps and then 
study how these traps contribute to the clock process. 

The coarse-graining procedure of~\cite{BCM06} studies the process $Y$ 
(resp.~$X$) only before the exit from a large disk $\mathbb D(n)$, 
$n\in \mathbb N$, of area $m 2^n n^{1-\alpha }$ centred at the 
origin. We denote by $R(n)$ the radius of this disk, 
$R_n=\sqrt{\pi^{-1}2^nn^{1-\alpha }}$. The random walk $Y$ makes 
$O(R_n^2)$ steps in $\mathbb D(n)$ and visits $O(R_n^2/\log R_n^2)$ 
different traps, as is well known. Therefore, the deepest traps that it visits 
have a depth of order 
$O((R_n^2/\log R_n^2)^{1/\alpha })=O(2^{n/\alpha }n^{-1})$. We 
therefore approximate the clock process by the time that $X$ spends 
in the  set
\begin{equation}
  T_\varepsilon^M=\{x\in \mathbb D(n):
    \varepsilon 2^{n/\alpha }n^{-1}\le \tau_x < M 2^{n/\alpha }n^{-1}\},
\end{equation}
where $\varepsilon $ is a small and $M$ a large constant. We call the 
traps in this set the \textit{deep traps}. It must be proved that 
this approximation is correct. We do not want, however, to deal with 
this problem here (see \cite{BCM06}).

After giving the proper depth scale, we can study how the deep traps are 
visited. This is where the coarse-graining is crucial. 
We cut the trajectory of the process $Y$ into short pieces. Every 
such piece of the trajectory
ends when $Y$ exits for the first time the disk of area 
$2^nn^\gamma $ around its initial point. At this moment a 
new piece starts.  Clearly, we should take $\gamma < 1-\alpha $. 
Formally, we set $j_0^n=0$, and then we define 
recursively 
\begin{equation}
  j_i^n=\min\big\{k>j_{i-1}^n:\dist(Y(k),Y(j^n_{i-1}))\ge
    \sqrt{\pi^{-1}2^nn^\gamma} \big\}.
\end{equation}
We use $x_i^n$ to denote the starting points of the pieces of the 
trajectory, $x^n_i=Y(j_i^n)$. It can be seen easily that the 
number of pieces of the trajectory  before the exit from 
$\mathbb D(n)$ is of order $O(n^{1-\alpha -\gamma })$. 

We look at the time that the walk spends in the deep traps during one 
piece of trajectory: we define the \textit{score} of the piece $i$ by
\begin{equation}
  s_i^n=\sum_{k=j^n_{i}}^{j^n_{i+1}-1}e_k\tau_{Y(k)}
  \bbone\{Y(k)\in 
    T_\varepsilon^M\}.
\end{equation}

To study the scores it is convenient to introduce another family of 
random variables $s_x^n$ indexed by $x\in \mathbb D(n)$. We set the 
distribution of $s_x^n$ to be the same as the distribution of $s^n_i$ 
conditioned on the fact that the $i^{\text{th}}$ piece of the trajectory 
starts at $x$, i.e.~conditioned on $x_i^n=x$. The main technical piece of the proof 
is to show that the law of $s_x^n$ does not depend on $x$ (with a small 
  exceptional set): 
\begin{lemma}
  \label{l:coarse}
  Fix $\kappa >0$ large enough  and define
  \begin{equation}
    \mathcal  E(n)=\{x\in \mathbb D(n): \dist(x,T_\varepsilon^M)\ge \sqrt{\pi^{-1}2^nn^{-\kappa }}\}. 
  \end{equation}
  Then for $\mathbb P$-a.e.\ random environment $\boldsymbol \tau $,
  uniformly for $x\in \mathcal E(n)$
  \begin{equation}
    \label{e:minmax}
    \lim_{n\to\infty}
    \frac{1-\mathbb E[\exp(- 
	  \frac{\lambda s_x^n}{2^{n/\alpha }})
	|s_x<\infty,\boldsymbol \tau ]}
    {n^{\alpha +\gamma -1}}
    =F(\lambda ). 
  \end{equation}
  Here
  \begin{equation}
    F(\lambda )=F(\lambda ;\varepsilon ,M,\alpha )=
    \mathcal  K\bigg(p_\varepsilon^M -
      \int_\varepsilon^M
      \frac{\alpha }
      {1+\mathcal K' \lambda z }\cdot
      \frac{1}{z^{\alpha +1}}\, \d z\bigg)
  \end{equation}
  with $\mathcal K'=\pi^{-1}\log 2$ and $\mathcal K=(\log 2)^{-1}$.
\end{lemma}
This lemma is a consequence of the following four facts, whose proofs 
are based on classical sharp estimates on the Green's function for 
the simple random walk on $\mathbb Z^d$ and on certain ``homogeneity'' 
properties of the random environment. These proofs can be found in 
\cite{BCM06}.

1. $s_x^n$ is equal to $0$ with probability 
$1-\mathcal K p_\varepsilon^M n^{\alpha +\gamma -1}(1+o(1))$, where 
$p_\varepsilon^M=\varepsilon^{-\alpha }-M^{-\alpha }$. That means 
that typically no deep trap is visited in a piece of the trajectory. 
Further, it implies that only a finite number of pieces has a non-zero score  
before the exit from $\mathbb D(n)$. 

2. With probability 
$\mathcal K p_\varepsilon^M n^{\alpha +\gamma -1}(1+o(1))$ the random 
walk visits (many times) only one deep trap during one piece, call it 
$y$. The probability that two or more deep traps are visited during 
one piece is $O(n^{2(\alpha +\gamma -1)})$, therefore, with 
overwhelming probability,  this event does not occur before the exit 
from~ $\mathbb D(n)$. 

3. In the case when one deep trap $y$  is visited, the distribution of its  
normalised depth  $2^{-n/\alpha } n \tau_y$ converges to the 
distribution on $[\varepsilon ,M]$ with the density $p(u)$ 
proportional to~$u^{-\alpha -1}$. 

4. The number of visits to $y$  is geometrically distributed with mean 
$\mathcal K' n(1+o(1))$. Therefore, conditionally on hitting $y$ the 
score has an exponential distribution with mean 
$\mathcal K' n \tau_y$ which is of order $O(2^{n/\alpha })$.

\medskip

Using Lemma~\ref{l:coarse}, it is not difficult to check that the 
scores $s_i^n$ are asymptotically i.i.d. Actually, to transfer the 
uniform convergence of distributions of $s_x^n$ into the asymptotic 
i.i.d.~property of $s_i^n$ it is sufficient to check that   with an 
overwhelming probability all pieces of trajectory before the exit of 
$\mathbb D(n)$ do start in $\mathcal E(n)$. 

The asymptotic i.i.d.~property then yields the convergence of 
the normalised sum of scores, 
$2^{-n/\alpha }\sum_{i=0}^{t n^{1-\alpha -\gamma }}s_i^n$, to a L\'evy 
process. The L\'evy measure of this process can be computed from 
Lemma~\ref{l:coarse}. Using the knowledge of the L\'evy measure, one can then
prove that for any $T>0$ it is possible to choose $m$ large enough 
such that as $\varepsilon \to 0$ and $M\to \infty$ the distribution 
of this L\'evy process on $[0,T]$ approaches the distribution  of  an 
$\alpha $-stable subordinator on $[0,T]$. The convergence of the clock 
process to the same subordinator then follows since the sum of scores 
is a good approximation of the clock process.

\subsection{References} 

The trap model in dimension $d\ge 2$ is studied in \cite{MB96,RMB01}. 
In the case $a=0$, a mathematical proof of aging has been given in 
\cite{Cer03,BCM06}. This proof is based on the  coarse-graining of 
the trajectories of the BTM sketched in the last section. In 
\cite{BC06} we establish the fractional-kinetics scaling limit based 
on the arguments of \cite{BCM06}. The case where $a\neq 0$ is 
discussed in \cite{RMB01}, but is  still an open problem. 

Let us mention that the clock process introduced in 
Definition~\ref{d:clock} is close to the problem of Random Walk in 
Random Scenery (RWRS), except for the extra randomisation due to 
exponential waiting times. However, here the tails of the scenery 
distribution are heavier than in the recent works on RWRS. 

In dimension $d=1$, RWRS with heavy tails have been studied by Kesten 
and Spitzer \cite{KS79}. The F.I.N.~diffusion could  indeed be seen 
as a Brownian motion time-changed by a (dependent) Kesten-Spitzer 
clock. 

\section{The arcsine law as a universal aging scheme} 
\label{s:general}
In this section we explain a general strategy that can be used to 
prove aging of the functions $\Pi $ and $R$ in the BTM on many 
different graphs (including $\mathbb Z^d$ for $d\ge 2$,   tori in 
  $\mathbb Z^d$,   large complete graphs, and high-dimensional 
  hypercubes). When this strategy can be used, the 
behaviour  of $R(t_w,t_w+\theta t_w;\boldsymbol \tau )$ for large 
$t_w$ can be expressed using the distribution function of the 
generalised arcsine 
distribution with parameter~$\alpha $, $\Asl_\alpha (\cdot)$. 

More formally, we will consider in this section a sequence of 
the Bouchaud trap models $\{\BTM(G_n,\boldsymbol \tau_n,0):n\in \mathbb N\}$. 
We will prove that for properly chosen time scales $t_w(n)$ 
\begin{equation}
  \lim_{n\to\infty}R_n(t_w(n),(1+\theta )t_w(n);\boldsymbol \tau 
    )=\Asl_\alpha (1/1+\theta ),
\end{equation}
where 
\begin{equation}
  \label{e:afdaf}
    R_n(t_w,t_w+t;\boldsymbol \tau ):=\mathbb P[X_n(t_w)=X_n(t_w+t)|\boldsymbol \tau ].
\end{equation}
The possible time scales $t_w(n)$ will depend on the graphs  $G_n$ 
and the laws of the depths $\boldsymbol \tau_n$. We will always 
try to prove aging results on the widest possible range of time scales. 

We will consider only the case of RHT dynamics, i.e. $a=0$. 
So that, the embedded discrete-time process $Y_n$ is the simple 
random walk on $G_n$, and the continuous time Markov chain $X_n$ is a 
time change of $Y$.  Our strategy relies on an approximation of the clock 
process (see Definition~\ref{d:clock}),
 more precisely of its rescaling,  by an $\alpha $-stable 
subordinator.  We will then show that the event 
$X_n(t_w)=X_n((1+\theta )t_w)$ is approximated by the event that 
the subordinator jumps over the interval $[1,1+\theta ]$. The 
classical arcsine law for L\'evy processes (see 
  Proposition~\ref{p:stablejump} in the Appendix) will then imply the 
aging result~\eqref{e:afdaf}.

\subsection{Aging on  large complete graphs}
\label{ss:complete}
As a warm-up, and in order to explain our strategy, we give here a 
complete proof of  aging on a large complete graph in this 
sub-section. The advantage of this graph is that the embedded simple 
random walk is particularly simple. On the other hand, almost all 
effects of the trapping landscape are already present. The method of 
the proof that we use here is probably not the simplest one. Its main 
ideas can however be adapted to more complex graphs; the structure of 
the proof stays the same, but a relatively fine control of 
the simple random walk on these graphs is then required.

Note that the complete graph is the graph that was proposed  in the 
original paper of Bouchaud \cite{Bou92}. The aging on this graph is 
proved using renewal arguments in \cite{BD95}, see also the introduction 
to \cite{BBG03b}. Another proof of aging and much more can be found 
in \cite{BF05} where eigenvalues and eigenvectors of the generator of 
$X$ are very carefully analysed. 

Let $G_n$ be a complete graph with $n$ vertices, 
$G_n = (\mathcal V_n, \mathcal E_n)$, where 
$\mathcal V_n=\{1,\dots ,n\}$ and 
$\mathcal E_n=\<x,y\>, x,y\in \mathcal V_n$. Note that we  include 
loops $\<x,x\>$ into the graph, so that jumps from $x$ to $x$ are 
possible. This makes the embedded random walk of the RHT dynamics 
extremely simple: the positions $Y_n(i)$, $i\in \mathbb N$, are 
i.i.d.~uniform random variables on $\mathcal V_n$. We also suppose that the starting 
position of the process, $X_n(0)=Y_n(0)$, is uniformly distributed on 
$\mathcal V_n$. Finally, we assume that the trapping landscape is 
given by an i.i.d.~sequence $\{\tau_x:x\in \mathbb N\}$ independent 
of $n$ such that Assumption~\ref{a:limit} holds.

As the process $X_n$ is a time change of the i.i.d.~sequence $Y_n(i)$, 
we need to study the clock process $S_n(t)$. Observe that  
  $S_n(k)$ for  $k\in \mathbb N$ is the time of $k^{\text{th}}$ jump of $X_n$, 
\begin{equation}
  S_n(k)=\sum_{i=0}^{k-1}\tau_{Y_n(i)}e_i,
\end{equation}
where $\{e_i\}$ is an i.i.d.~sequence of exponential random variables 
with mean one.  

Since we have included the loops into the graph, we redefine the two-time 
function $\Pi $ slightly:
\begin{equation}
  \Pi_n(t_w,t_w+t;\boldsymbol \tau ) := \mathbb P
  \big[\{S_n(k):k\in\mathbb  N\}\cap [t_w,t_w+t]=\emptyset\big|\boldsymbol \tau \big].
\end{equation}
Note that this definition differs from the original one only if a 
jump from $X(t_w)$ to $X(t_w)$ occurs. The probability of this event is 
$1/n$ and is thus negligible for large $n$. Similarly, for the 
complete graph the function $R_n$ differs very little from the 
function $\Pi_n$, a correction is again of  order $1/n$. 

We prove aging on time scales  smaller than $n^{1/\alpha }$. 
The scale $n^{1/\alpha }$ appears  because the deepest trap in 
$\mathcal V_n$ has a depth of this order as can be easily verified 
from Assumption~\ref{a:limit}. Therefore, at time scales shorter than 
$n^{1/\alpha }$ the process has not enough time to reach the 
equilibrium, and aging can be observed: 

\begin{theorem}[Quenched aging on the complete graph]
  \label{t:agingcg}
  Let  $0<\kappa < 1/\alpha $ and let  $t_w=t_w(n)=n^{\kappa }$.  Then 
  for a.e.~random environment $\boldsymbol \tau $
  \begin{equation}
    \lim_{n\to\infty}
    \Pi_n(t_w(n),(1+\theta )t_w(n);\boldsymbol \tau )=
    \Asl_\alpha (1/1+\theta ).
  \end{equation}
  Moreover, the clock process converges to the stable subordinator $V$
  with the  L\'evy measure $\alpha \Gamma (1+\alpha )u^{-\alpha -1}\d u$, 
  \begin{equation}
    \frac {S_n(sn^\kappa )}{n^{\kappa /\alpha }}\xrightarrow{n\to\infty}
    V(s).
  \end{equation}
\end{theorem}

\begin{remark}
  1. A similar result holds on the shortest possible time scale 
  $\kappa =0$: if $t_w(n)=t n^0 = t$, then 
  \begin{equation}
    \lim_{t\to\infty}\lim_{n\to \infty}
    \Pi_n(t_w,(1+\theta )t_w;\boldsymbol \tau )=
    \Asl_\alpha (1/1+\theta ).
  \end{equation}
  This is actually proved in \cite{BD95} and \cite{BF05}.

  2. For the longest possible time scales 
  $t_w(n)=t n^{1/\alpha }$ a double limiting procedure is also necessary. 
  Essentially the same arguments 
  that are used to prove Theorem~\ref{t:agingcg} yield that 
  $\Pi_n(t_w(n),(1+\theta )t_w(n);\boldsymbol \tau )$ converges to 
  $\Asl_\alpha (1/1+\theta )$ in probability, that is for any 
  $\varepsilon >0$
  \begin{equation}
    \label{e:cginproba}
    \lim_{t\to 0}\lim_{n\to\infty}
    \mathbb P\Big[
    \Big|\Pi (tn^{1/\alpha },(1+\theta )tn^{1/\alpha };\boldsymbol \tau 
      )-\Asl_\alpha (1/1+\theta )\Big|\ge \varepsilon \Big]=0.
  \end{equation}

  3. No scaling limit for $X$ exists since there is no such thing as 
  a scaling limit for the complete graph.
\end{remark}

To prove Theorem~\ref{t:agingcg} we consider a rescaling of the time 
change process, 
\begin{equation}
  \label{e:lpdef}
  \mathcal S_n(t)= \frac 1 {n^\kappa } S_n(k) \qquad
  \text{for $t\in [u_n(k),u_n(k+1))$,}
\end{equation}
where $u_n(k)=n^{-\kappa \alpha }\sum_{i=1}^k e'_i$, and 
$\{e'_i,i\in \mathbb N\}$ is another sequence of mean one 
i.i.d.~exponential random variables. We introduce this sequence in 
order to randomise the times of jumps of $\mathcal S_n$. Thanks to 
this randomisation $\mathcal S_n$ is the L\'evy process. The standard 
rescaling of $S_n$, $n^{-\kappa } S_n( t n^{\kappa \alpha })$ does 
not have this property. $\mathcal S_n$ is compound Poisson process. 
The intensity of its jumps is $n^{\alpha \kappa }$. Every jump has 
the same distribution as $n^{-\kappa }\tau_x e_i$, where $x$ is 
uniformly distributed in $\mathcal V_n$. Therefore,  
the L\'evy measure $\mu_n$ of $\mathcal S_n$ is 
\begin{equation}
  \mu_n(\d u)= n^{\alpha \kappa }
  \frac 1n \sum_{x\in \mathcal V_n}
  \frac {n^{\kappa }e^{-un^\kappa /\tau_x}}{\tau_x} \,\d u =
  n^{\alpha \kappa +\kappa -1}\!\sum_{x\in \mathcal V_n}
  \frac {e^{-u n^\kappa /\tau_x}} {\tau_x} \,\d u.
\end{equation}
It follows from definitions of $\Pi_n$ and $\mathcal S_n$ that 
\begin{equation}
  \label{e:hhj}
  \Pi_n (n^\kappa ,(1+\theta )n^\kappa ;\boldsymbol \tau )=
  \mathbb P\big[\{\mathcal S_n(t):t\in \mathbb R\}
    \cap [1,1+\theta ]=\emptyset\big],
\end{equation}
that is the probability that $\mathcal S_n$ jumps over the interval 
$[1,1+\theta ]$. 

The idea of the proof is the following. We write $\mathcal S_n$ as a 
sum of three independent L\'evy processes, $\mathcal S_{n,M}=\mathcal S_{n,M}^\infty$, 
$\mathcal S_{n,\varepsilon}^M$ and 
$\mathcal S_n^\varepsilon=\mathcal S_{n,0}^\varepsilon  $ with the 
L\'evy measures $\mu _{n,M}=\mu_{n,M}^\infty$, 
$\mu _{n,\varepsilon}^M$ and 
$\mu_n ^\varepsilon=\mu _{n,0}^\varepsilon$, where
\begin{equation}
  \mu_{n,a}^b(\d u)=
  n^{\alpha \kappa +\kappa-1}\sum_{x\in \mathcal V_n}
  \bbone\{x\in T_a^b(n)\}
  \frac 1 {\tau_x} e^{-u n^\kappa /\tau_x}\,\d u, 
\end{equation}
and
\begin{equation}
  T_a^b(n)=\Big\{x\in \mathcal V_n:
    \frac {\tau_x}{n^\kappa }\in [a,b)\Big\}
\end{equation}
That is we divide traps into three categories: (a) \textit{the very 
  deep traps}, $x\in T_M(n)$, i.e.~$\tau_x\ge Mn^\kappa $, (b) 
\textit{the deep traps}, $x\in T_\varepsilon^M(n)$, 
i.e.~$\tau_x/n^\kappa\in [\varepsilon ,M)$, (c) \textit{the shallow 
  traps}, $x\in T^\varepsilon $, i.e.~$\tau_x<\varepsilon n^\kappa$. 
We consider the contributions of these different categories to the 
clock process separately.

We show that  the contribution of the deep traps, i.e.~the process 
$\mathcal S_{n,\varepsilon}^M$, is well approximated by an $\alpha $-stable 
subordinator, at least if $\varepsilon $ is small and $M$ large 
enough. The proof of this fact relies on the weak convergence of the 
L\'evy measures  $\mu^M_{n,\varepsilon} $.  
We prove this convergence only on the interval $[0,T]$ where $T$ is 
chosen such that $\mathcal S_{n,\varepsilon}^M$ is larger 
than $(1+\theta )$ with a large probability. 

Further we prove that the very deep and shallow traps 
can be almost neglected. More exactly, we prove that $\mathcal S_{n,M}(T)=0$ 
with a large probability, that is $\mathcal S_{n,M}$ does not jump before 
$T$, and that $\mathcal S_n^\varepsilon (T)$ can be made small by choosing 
$\varepsilon $ small enough. These fact together will imply that 
$\mathcal S_n$ is well approximated on $[0,T]$ by an $\alpha $-stable 
subordinator, and the claim of the theorem  follows from 
\eqref{e:hhj} and the arcsine law for stable subordinators, more 
precisely from Corollary~\ref{c:stablejump}. 

\subsubsection{Deep traps}

The following proposition describes the contribution of the deep traps to 
the time change, that is the process $\mathcal S_{n,\varepsilon}^M$. 

\begin{proposition}
  \label{p:afdas}
  The L\'evy measures $\mu_{n,\varepsilon }^M$ of $\mathcal S_{n,\varepsilon}^M$ 
  converge weakly as $n\to \infty$ to the measure $\mu_\varepsilon^M$ given by 
  \begin{equation}
    \label{e:integral}
    \mu_\varepsilon^M(\d u)=
    \int_\varepsilon^M \frac \alpha  {z^{\alpha+2}}\, e^{-u/z} \,\d z 
    \, \d u.
  \end{equation}
\end{proposition}

This proposition has an important corollary that states that 
$\mathcal S_{n,\varepsilon}^M$ can be well approximated for large  $n$ by 
an $\alpha $-stable subordinator. To this end, let 
$\mathcal Z_\varepsilon^M$ be a subordinator with the L\'evy measure
\begin{equation}
  \tilde \mu_\varepsilon^M(\d u)=
  \int_0^\varepsilon  \frac \alpha  {z^{\alpha +2}}\, e^{-u/z} \,\d z 
    \, \d u
    + \int_M^\infty  \frac \alpha  {z^{\alpha +2}}\, e^{-u/z} \,\d z 
    \, \d u
\end{equation}
independent of all already introduced random variables. 
\begin{corollary}
  \label{c:stableapp}The process 
  $\mathcal S_{n,\varepsilon}^M + \mathcal Z_\varepsilon^M$ converge 
  weakly in the Skorokhod topology on $D([0,T],\mathbb R)$ to an 
  $\alpha $-stable subordinator. Moreover, for any $T$ and $\delta >0$ 
  it is possible to choose $\varepsilon $ small and $M$ large enough 
  such that 
  \begin{equation}
    \label{e:maly}
    \mathbb P[\mathcal Z_\varepsilon^M(T)>\delta ]<\delta  .
  \end{equation}
\end{corollary}
\begin{proof}[Proof of Corollary~\ref{c:stableapp}]
  The first claim is consequence  of the convergence of L\'evy measures 
  (Proposition~\ref{p:afdas}) 
  together with Lemma~\ref{l:subconv}. The fact that the limit is a stable
  subordinator follows from 
  \begin{equation}
    \tilde \mu_\varepsilon^M(\d u)+\mu_\varepsilon^M(\d u)= 
    \alpha \Gamma (1+\alpha ) u^{-\alpha -1}\,\d u. 
\end{equation}
  Finally, to prove \eqref{e:maly} observe that for all $\lambda >0$
  \begin{equation}
    \mathbb E[e^{-\lambda \mathcal Z_\varepsilon^M(T)}]=
    e^{-T \int_0^\infty (1-e^{-\lambda x})\tilde \mu_\varepsilon^M(\d x)}
    \xrightarrow{\varepsilon \to 0, M\to\infty } 1,
  \end{equation}
  that is the law of $\mathcal Z_\varepsilon^M(T)$ converges weakly to 
  the Dirac mass at $0$.
\end{proof}

To prove Proposition~\ref{p:afdas} we first show that for 
a.e.~$\boldsymbol \tau $ there is the right number 
of deep traps with depths approximately $un^\kappa $ in $\mathcal V_n$, 
$u\in [\varepsilon ,M)$. This follows from a ``law-of-large-number-type'' 
argument using the fact that we have a large number of traps 
with such depth. This property fails to be true if  
$\kappa \ge 1/\alpha $ (which explains why we cannot get a.s. result 
  for $t_w=cn^{1/\alpha }$ and why we need the double limit procedure).

\begin{lemma}
  \label{l:rightnumber}
  Let $0<\kappa < 1/\alpha $. Then there exists a function 
  $h(n)$, $h(n)\to 0$ as $n\to \infty$,  
  such that 
  \begin{equation}
    \lim_{n\to \infty}
    \sup_{u\in [\varepsilon ,M)\cap \{\mathbb Z/h(n)\}}
    \bigg|
    \frac {|T_u^{u+h(n)}(n)|}{h(n) n^{1-\alpha \kappa }}-
    \alpha u^{-\alpha -1}\bigg|=0,\qquad 
    \text{$\boldsymbol \tau$-a.s.}
  \end{equation}
\end{lemma}

\begin{proof}
  Let $g(u)$ be defined by 
  $\mathbb P[\tau_x\ge u]=u^{-\alpha }(1+g (u))$. It follows from 
  Assumption~\ref{a:limit} that  $\lim_{u\to \infty}g(u)=0$. Take 
  $h(n)$ such that $h(n)\to 0$, $h(n)\ge (\log n)^{-1}$, 
  $h(n)\gg \sup_{u\in [\varepsilon ,M)} g(un^\kappa )$. Then, as
  $n\to\infty$,
  \begin{equation}
    \label{e:taupr}
    \mathbb P\big[\tau_x\in [un^\kappa ,(u+h(n))n^\kappa )\big]=\alpha 
    u^{-1-\alpha }h(n) n^{-\alpha \kappa }(1+o(1)).
  \end{equation}
  Indeed, by definition of $g$, the left-hand side is equal to 
  \begin{equation}
    u^{-\alpha }n^{-\kappa \alpha } 
    \Big[1-\Big(1+\frac {h(n)}{u}\Big)^{-\alpha }\Big]+
    u^{-\alpha }n^{-\kappa \alpha }
    \Big[g(un^\kappa )-
      g((u+h(n))n^\kappa )\Big(1+\frac {h(n)}u\Big)^{-\alpha }\Big].
  \end{equation}
  The last expression is equal to the right-hand side of 
  \eqref{e:taupr} as follows from the definitions of $g$ and $h$.

  Fix $u\in[\varepsilon ,M]$ and $\varepsilon >0$. Using twice the exponential 
  Chebyshev inequality together with \eqref{e:taupr} we get 
  \begin{equation}
    \begin{aligned}
      \mathbb P\Big[\frac{|T_u^{u+h(n)}(n)|}{h(n) n^{1-\alpha \kappa }}
	\notin
	(1-\varepsilon ,1+\varepsilon )\alpha u^{-\alpha -1}\Big]
      \le
      C \exp\{-c(\varepsilon ) u^{-\alpha -1} n^{1-\alpha \kappa }h(n)\}.
  \end{aligned}
  \end{equation}
  Summing over $u\in [\varepsilon ,M)\cap \{\mathbb Z/h(n)\}$, the 
  claim of the lemma then follows using the 
  Borel-Cantelli Lemma and the fact that $h(n)\ge (\log n)^{-1}$.
\end{proof}

\begin{proof}[Proof of Proposition~\ref{p:afdas}]
  By definition of $\mu_{n,\varepsilon }^M$,
  \begin{equation}
    \begin{aligned}
      \mu_{n,\varepsilon }^M(\d u)&=
      n^{\alpha \kappa +\kappa -1}\sum_{x\in  T_\varepsilon^M(n)}
      \frac 1 {\tau_x} e^{-u n^\kappa /\tau_x}\,\d u
      \\&=
      n^{\alpha \kappa +\kappa -1}
      \sum_{i=0}^{(M-\varepsilon )/h(n)}
      \sum_{x\in T_{ih(n)}^{(i+1)h(n)}}
      \frac {e^{-u /ih(n)}} {n^\kappa i h(n)} \,\d u(1+o(1))
      \\&=
      n^{\alpha \kappa  -1}
      \sum_{i=0}^{(M-\varepsilon )/h(n)}
      |T_{ih(n)}^{(i+1)h(n)}(n)|
      \frac 1 { i h(n)} e^{-u /ih(n)}\,\d u(1+o(1)).
    \end{aligned}
  \end{equation}
  By the previous lemma, for large $n$,  $\boldsymbol \tau $-a.s 
  \begin{equation}
      \mu_{n,\varepsilon }^M(\d u)
      =
      \sum_{i=0}^{(M-\varepsilon )/h(n)}
      \alpha (i h(n))^{-\alpha -1}
      \frac {e^{-u /ih(n)}} {i h(n)} h(n)\,\d u(1+o(1)),
  \end{equation}
  which is the Riemann sum of the integral in \eqref{e:integral} with 
  the mesh size $h(n)$. Since $h(n)\to 0$ as $n\to \infty$, the proof is finished.
\end{proof}

\subsubsection{Shallow traps}
\label{sss:shallow}
We prove that the contribution of the shallow traps at any final 
instant $T$, 
$\mathcal S^\varepsilon_n(T)$, is small if $\varepsilon $ is chosen 
small enough. 

\begin{lemma}
  \label{l:cgshallow}
  There exists a large constant $K$ such that for all $T>0$ and $n$ large
  \begin{equation}
    \mathbb E[\mathcal S^\varepsilon_n(T)|\boldsymbol \tau ]
    <K T 
    \varepsilon^{1-\alpha },\qquad \boldsymbol \tau\text{-a.s.}
  \end{equation}
\end{lemma}
\begin{proof}
  Let $\xi $ be such that $T\in [u_n(\xi ),u_n(\xi +1))$ (see \eqref{e:lpdef}). It 
  is easy to see that $\xi $ has Poisson distribution, $\mathbb E[\xi ]=Tn^{\kappa \alpha }$. By 
  definition of $\mathcal S^\varepsilon_n$,
  \begin{equation}
    \mathcal S^\varepsilon_n(T)=\frac 1 {n^\kappa }
    \sum_{i=0}^\xi  e_i \tau_{Y_n(i)}
    \bbone \{Y_n(i)\in T^\varepsilon (n)\}.
  \end{equation}
  To bound its expected value we divide $T^\varepsilon (n)$ into 
  slices $T_{\varepsilon 2^{-i}}^{\varepsilon{2^{-i+1}}}$,
  $i\in \mathbb N$. We define
  \begin{equation}
    A_n(i):=\Big\{\boldsymbol \tau :
      \mathbb E \Big[\frac 1 {n^\kappa }
    \sum_{i=0}^\xi  e_i \tau_{Y_n(i)}
    \bbone \big\{Y_n(i)\in T_{\varepsilon{2^{-i}}}^{\varepsilon{2^{-i+1}}} (n)\big\}
    \Big| \boldsymbol \tau \Big]
  < K' \varepsilon^{1-\alpha }2^{i(\alpha -1)}\Big\}.
  \end{equation}
  We  show that there is $K'$ such that a.s.
  \begin{equation}
    \label{e:cgshallow}
    A_n:=\bigcap_{i=0}^\infty A_n(i)\quad 
    \text{occurs for all but finitely many $n$}.
  \end{equation}
  Lemma~\ref{l:cgshallow} is then a direct consequence of this claim. 

  To show \eqref{e:cgshallow} we first estimate the probability that 
  a fixed site, say 1, is in $T_{\varepsilon{2^{-i}}}^{\varepsilon{2^{-i+1}}}$,
  \begin{equation}
    \label{e:pni}
    p_{n,i}:=\mathbb P[1\in 
      T_{\varepsilon{2^{-i}}}^{\varepsilon{2^{-i+1}}}]\le 
    C\varepsilon^{-\alpha }2^{i\alpha } n^{-\kappa \alpha }
  \end{equation}
  as follows from Assumption~\ref{a:limit}. For any $x\in \mathcal V_n$,
  \begin{equation}
    \label{e:greencga}
     \mathbb E\Big[\sum_{i=0}^\xi  e_i \bbone\{Y_n(i)=x\}\Big]=
    Tn^{\kappa \alpha -1 }, 
  \end{equation}
  therefore
  \begin{multline}
      \mathbb P\Big[
	\mathbb E \Big[\frac 1 {n^\kappa }
	  \sum_{i=0}^\xi  e_i \tau_{Y_n(i)}
	  \bbone \big\{Y_n(i)\in 
	    T_{\varepsilon{2^{-i}}}^{\varepsilon{2^{-i+1}}} (n)\big\}
	  \Big| \boldsymbol \tau \Big]
	\ge K' \varepsilon^{1-\alpha }2^{i(\alpha -1)}\Big]
      \\\le
      \mathbb P\Big[ \sum_{x\in \mathcal V_n}
	T n^{\kappa \alpha -1}  \varepsilon  2^{-i+1}
	\bbone \big\{x\in 
	  T_{\varepsilon{2^{-i}}}^{\varepsilon{2^{-i+1}}} (n)\big\}
	\ge K' \varepsilon^{1-\alpha }2^{i(\alpha -1)}\Big].
  \end{multline}
  By Chebyshev inequality this is bounded from above by
  \begin{equation}
    \exp\{-\lambda n^{1-\kappa \alpha}   
      c K' \varepsilon^{-\alpha }2^{i\alpha }\} 
    \prod_{x\in \mathcal V_n}(p_{n,i}(e^\lambda -1)+1).
  \end{equation}
  Using $1+x\le e^x$ and \eqref{e:pni}, it is easy to see that for 
  $K'$ large enough this is smaller than
  $ \exp\{- c n^{1-\kappa \alpha } 2^{i\alpha } \varepsilon^{-\alpha }\}$.
  Therefore, summing over $i\in \mathbb N$,
  \begin{equation}
    \mathbb P[A_n^c]\le \exp\{-c' n^{1-\kappa \alpha } 
      \varepsilon^{-\alpha }\}
  \end{equation}
  An application of the Borel-Cantelli lemma finishes the proof of 
  \eqref{e:cgshallow} and therefore of Lemma~\ref{l:cgshallow}.
\end{proof}

\subsubsection{Very deep traps}
We show that, with  a large 
probability, none of the very deep traps is visited during the first 
$Tn^{\kappa \alpha }$ steps: 

\begin{lemma}
  \label{l:cgdeep}
  If $0<\kappa < 1/\alpha $, then for all $n$ large enough, $\boldsymbol \tau $-a.s.
  \begin{equation}
    \mathbb P[\mathcal S_{n,M}(T)\neq 0|\boldsymbol \tau ]\le c TM^{-\alpha }.
  \end{equation}
\end{lemma}
\begin{proof}
  Define $H=\inf\{k\in \mathbb N_0: Y_n(k)\in T_{M}(n)\}$. By 
  definition of $\mathcal S_{n,M}$, the claim of the lemma is 
  equivalent to 
  $\mathbb P[H\le \xi|\boldsymbol \tau ]\le c TM^{-\alpha }$, with the 
  same $\xi$ as in the previous proof. Since
  \begin{equation}
    \mathbb P[H\le \xi|\boldsymbol \tau ]\le 
    \mathbb P[H<2Tn^{\kappa \alpha }|\boldsymbol \tau ]+
    \mathbb P[\xi>2Tn^{\kappa \alpha }],
  \end{equation}
  and the second term on the right-hand side decreases exponentially 
  in $n$, it is sufficient to bound the first term. Using the fact 
  that $\boldsymbol \tau $-a.s.~$|T_M(n)|\le C M^{-\alpha }n^{1-\kappa \alpha }$ 
  (which can be proved 
  similarly as Lemma~\ref{l:rightnumber}), we have
  \begin{equation}
      \mathbb P[H<2Tn^{\kappa \alpha }|\boldsymbol \tau ]
      =
      \mathbb P\Big[\bigcup_{k=0}^{2Tn^{\kappa \alpha }}\{Y_n(k)\in 
	  T^M(n)\}\Big|\boldsymbol \tau \Big]
      \le
      2 Tn^{\kappa \alpha }\frac{ |T^M(n)|}n
      \le c M^{-\alpha } T.
  \end{equation}
  This finishes the proof.
\end{proof}

\subsubsection{Proof of Theorem~\ref{t:agingcg}}
\label{sss:agingcg}
We can now finish the proof of aging on a large complete graph. All 
claims here are valid $\boldsymbol \tau $-a.s.

First, for $\delta >0$ we fix $T>0$ such that (for all $n$ large enough)
\begin{equation}
  \mathbb P[\mathcal S_{n,\varepsilon }^M(T)+\mathcal Z_\varepsilon^M(T)
    \le 2+\theta|\boldsymbol \tau ]<\frac \delta 4,
\end{equation}
which is possible due to Corollary~\ref{c:stableapp}. Further, we 
use this corollary, Lemmas~\ref{l:cgshallow} and~\ref{l:cgdeep} to 
fix $\varepsilon $ and $M$ such that
\begin{equation}
  \mathbb P[\mathcal S_{n,M }(T)\neq 0|\boldsymbol \tau ]
  \le \frac \delta 4,\qquad
  \mathbb P\Big[\mathcal S^\varepsilon_n\ge \frac \delta 2 \Big
    |\boldsymbol \tau \Big]\le 
  \frac \delta 4
  \quad\text{and}\quad
    \mathbb P\Big[\mathcal Z_\varepsilon^M\ge \frac \delta 2\Big]\le
  \frac \delta 4.
\end{equation}
Let $B(n)$ be the intersection of all events from the two previous 
displays. It follows that  $\mathbb P[B(n)^c]\le \delta $.  Let $E(n)$ be the event whose 
probability we are trying to estimate (see~\eqref{e:hhj}):
\begin{equation}
 E(n):=\big\{ 
  \{\mathcal S_n(t):t\in \mathbb R\}
  \cap [1,1+\theta ]=\emptyset\big\}.
\end{equation}
On $B(n)$ we can approximate $\mathcal S_n(t)$ by 
$\mathcal S_{n,\varepsilon }^M(t)+\mathcal Z_\varepsilon^M(t)$:
\begin{equation}
  \label{e:vzd}
  B\implies \sup_{t\in[0,T]}|\mathcal S_n(t)-
  \mathcal S_{n,\varepsilon }^M(t)+\mathcal 
  Z_\varepsilon^M(t)|<\delta.
\end{equation}
Further, let 
$\mathcal R_{n,\varepsilon}^M=\{\mathcal S_{n,\varepsilon }^M(t)+\mathcal Z_\varepsilon^M(t):t\in \mathbb R\}$. 
We define
\begin{equation}
  \begin{aligned}
    G_1(n)&:=
    \big\{ \mathcal R_{n,\varepsilon }^M\cap
      \big((1-\delta ,1+\varepsilon )\cup(1+\theta -\delta ,1+\theta 
	  +\delta )\big)\neq \emptyset\},\\
    G_2(n)&:=
    G_1(n)^c \cap 
    \big\{ \mathcal R_{n,\varepsilon }^M
      \cap(1,1+\theta )
      \neq \emptyset\big\},\\
    G_3(n)&:=
    G_1(n)^c \cap 
    \big\{ \mathcal R_{n,\varepsilon }^M
      \cap(1,1+\theta )
      = \emptyset\big\}=G_1(n)^c\cap G_2(n)^c.
  \end{aligned}
\end{equation}
These three events form a partition of the probability space. The 
reason for this partition is the following. If $G_1(n)\cap B(n)$  
happens, then the process 
$\mathcal S_{n,\varepsilon }^M+\mathcal Z_\varepsilon^M$ intersects 
$\delta $-neighbourhoods of  $1$ or $1+\theta $, and it is therefore 
not possible to decide if $E(n)$ is true. On the other hand,
\begin{equation}
  B(n)\cap G_3(n)\implies E(n)\qquad \text{and}\qquad
  B(n)\cap G_2(n)\implies E(n)^c
\end{equation}
as can be seen from \eqref{e:vzd}.
Therefore
\begin{equation}
  \label{e:bounds}
  \mathbb P[B(n)\cap G_3(n)|\boldsymbol \tau ]\le 
  \mathbb P[E(n)|\boldsymbol \tau ]\le
  \mathbb P[B(n)^c|\boldsymbol \tau ]
  +\mathbb P[G_1(n)|\boldsymbol \tau ]
  +\mathbb P[G_3(n)|\boldsymbol \tau ].
\end{equation}
Thanks to Corollary~\ref{c:stableapp} it is possible to 
estimate the probabilities of $G_1(n)$ and $G_3(n)$. Indeed, 
$\mathcal S_{n,\varepsilon }^M+\mathcal Z_\varepsilon^M$ converge 
weakly in the Skorokhod topology to an $\alpha $-stable subordinator. The 
probability that the subordinator hits any of the boundary points $1\pm \delta $, 
$1+\theta \pm \delta $ is zero. Therefore, it follows from the weak 
convergence and Corollaries~\ref{c:stableapp} and~\ref{c:stablejump} that
\begin{equation}
  \lim_{n\to\infty} \mathbb P[G_3(n)|\boldsymbol \tau ]=\Asl_\alpha \Big(\frac{1-\delta 
    }{1+\theta +\delta }\Big).
\end{equation}
Similarly, by \eqref{e:aslass},
\begin{equation}
  \label{e:bounds2}
  \mathbb P[G_1(n)|\boldsymbol \tau ]\le 1-\Asl_\alpha \Big(\frac {1-\delta }{1+\delta }\Big)+1-
  \Asl_\alpha \Big(\frac{1+\theta -\delta }{1+\theta +\delta }\Big)\le
  C\delta^{1-\alpha }.
\end{equation}
Since $\delta $ is arbitrary, the claim of the theorem follows  from 
\eqref{e:bounds}--\eqref{e:bounds2} and from the continuity of 
$\Asl_\alpha(\cdot) $. This proves aging on the complete graph for 
$\kappa < 1/\alpha $.

\subsection{The $\alpha $-stable subordinator as a universal clock} 
\label{ss:universal}
We 
will now give a general set of conditions which ensures that the result and 
the proof that we have given in the simple case of the complete graph 
apply. Consider  an arbitrary sequence of graphs $G_n$ and the Bouchaud 
trap models on them, $\BTM(G_n,\boldsymbol \tau_n,0)$.  We set 
$\nu =1$, therefore at  
every vertex $x$ the Markov chain $X$ waits an exponentially 
distributed time with mean $\tau_x/d_x$ ($d_x$ is the degree of $x$) 
and then it jumps to one of the 
neighbouring vertices with an equal probability. We suppose that 
Assumption~\ref{a:limit} holds, in particular that means that $\tau_x$ 
are i.i.d.  

We want to prove aging for this model at a time scale $t_w=t_w(n)$ 
using the same strategy as for the complete graph. That 
means, depending on the sequence $G_n$ and the time scale $t_w$:
\begin{itemize}
  \item to divide traps into three groups: shallow, deep, and very deep,
  \item to prove that the shallow  traps can be ignored since the 
  time spent there is negligible,
  \item to prove that the very deep traps can be ignored because they 
  are not visited in the proper time scale,
  \item to show that the contribution of the deep traps to the time 
  change can be approximated by a stable subordinator, which will 
  show the convergence of the clock process,
  \item and finally, to deduce aging for the two-time function $R$ 
  from this convergence.
\end{itemize}
The idea is 
that,  as in the case of the complete graph, 
the time change $S_n(j)$ should be dominated by a relatively 
small number of large contributions coming from the deep traps.

On the other hand we want to stay as general as possible: we do not 
want to  use any particular properties of the graph. Therefore, we 
formulate six conditions on the sequence $G_n$ and the simple random 
walk $Y_n(\cdot)$ on it. If these conditions are verified, the proof 
of aging can be finished in the spirit of Section~\ref{sss:agingcg}.
Proving these conditions should be dependent on the graphs $G_n$.

\medskip

To formulate the conditions it is necessary to choose several 
objects that depend on the particular sequence $G_n$ and on the 
observation time scale $t_w(n)$. 

First, it is necessary to fix a (random) time $\xi_n$  up to which we 
observe $Y_n$.  Second, a scale $g(n)$ for deep traps should be 
chosen according to $G_n$ and $t_w(n)$. This scale defines the set of 
the deep traps by
\begin{equation}
  \label{e:deftop}
  T_\varepsilon^M(n):=\{x\in \mathcal V_n: \varepsilon g(n)\le 
    \tau_x < M g(n)\}.
\end{equation}
A possible generalisation of the definition~\eqref{e:deftop} is described 
in the remark after the Theorem~\ref{t:main}.

For the complete graph we used $g(n)=n^{\kappa }$ and 
$\xi_n\sim n^{\kappa \alpha }$ (see  the proof of Lemma~\ref{l:cgshallow}).

As before, we use $T_M(n)=\{x:\tau_x\ge Mg(n)\}$ to denote the set of {\em very 
  deep traps}. Similarly, we write 
$T^\varepsilon (n)=\{x:\tau_x<\varepsilon g(n)\}$ for the set of 
{\em shallow traps}.

It should be possible to almost ignore these two sets. This is 
ensured by the following two conditions. Compare them with 
Lemmas~\ref{l:cgshallow} and~\ref{l:cgdeep}.
\begin{condition}
  \label{c:shallow}
  There is a function $h(\varepsilon )$  satisfying 
  $\lim_{\varepsilon \to 0} h(\varepsilon )=0$, such that for 
  a.e.~realisation of $\boldsymbol \tau $  
  and for all $n$ large enough
  \begin{equation}
    \mathbb E\Big[\sum_{i=0}^{\xi_n}e_i \tau_{Y_n(i)}\bbone\{Y_n(i)\in 
	T^\varepsilon (n)\}\Big|\boldsymbol \tau \Big]\le h(\varepsilon )t_w(n).
  \end{equation}
  In words, the expected time spent in the shallow traps before $\xi_n$ is 
  small with respect to $t_w(n)$.
\end{condition}

Let $H_n(A)$, $A\subset \mathcal V_n$,  
denotes the hitting time of $A$ by the simple random walk $Y_n$,
\begin{equation}
  H_{n}(A):=\inf\{i\ge 0:Y_n(i)\in A\}.
\end{equation}
\begin{condition}
  \label{c:deep}
  Given $\xi_n$, for any $\delta >0$ there exists $M$ large enough such 
  that for a.e.~realisation of $\boldsymbol \tau $ and for all $n$ large 
  \begin{equation}
    \mathbb P\big[H_n(T_M(n))\le \xi_n\big|\boldsymbol \tau \big]\le \delta.
  \end{equation}
\end{condition}

We should now ensure that the contribution of the deep traps to the 
time change can be approximated by an $\alpha $-stable subordinator. The 
following facts were crucial for the proof in the complete graph 
case: asymptotically, as $n\to\infty$,
\begin{itemize}
  \item every time a deep trap is visited, its depth is independent of 
  depths of the previously visited deep traps.
  \item the probability that a deep trap with depth $u g(n)=u n^\kappa $ is 
  visited is proportional to $u^{-\alpha }$.
  \item for fixed $\varepsilon $ and $M$ only a finite number of deep 
  traps was visited before the time-horizon $\xi $.
\end{itemize}
However, these facts could not be true for a general graph. Here, 
after leaving a deep trap, the process $Y$ typically hits this 
trap with larger probability than any other deep trap. On recurrent 
graphs, such as $\mathbb Z^2$, it even visits this trap a number of 
times that diverges with $n$.  

To overcome this problem we ``group'' the successive visits of one 
deep trap. The time spend during these visits there will then be 
considered as one contribution to the clock process.
We define $r_n(j)$ as the sequence of times when 
a new deep trap is visited, $r_n(0)=0$, and 
\begin{equation}
  r_n(i)=\min\big\{j>r_n(i-1): Y_n(j)\in T_\varepsilon^M(n)
    \setminus\{Y_n(r_n(i-1))\}\big\}. 
\end{equation}
We use $\zeta_n$ to denote the largest $j$ such that $r_n(j)\le \xi_n$,
\begin{equation}
  \label{e:defzeta}
  \zeta_n:=\max \{j:r_n(j)\le \xi_n\}.
\end{equation} 
We define the process $U_n(j)$ that records the trajectory of $Y_n$ (and 
  thus of~$X_n$) restricted to the deep traps, 
\begin{equation}
  \label{e:Undef}
  U_n(j):=Y_n(r_n(j)),\qquad j\in \mathbb N_0.
\end{equation}
Finally, define the score $s_n(j)$ be the time that $X_n$ spends at 
site $U_n(j)$ between steps $r_n(j)$ and $r_n(j+1)$, 
\begin{equation}
  \label{e:defsn}
  s_n(j):=\sum_{i=r_n(j)}^{r_n(j+1)}e_i \tau_{Y_n(i)} 
  \bbone\{Y_n(i)=U_n(j)\},\qquad
  \text{$j<\zeta_n$}.
\end{equation}
Denoting by $G_A^n(x,y)$ the Green's function of the simple random walk 
$Y_n$ killed on the first visit to a set $A\subset \mathcal V_n$, it 
is easy to observe that $s_n(j)$ has the exponential distribution 
with mean 
\begin{equation}
  \label{e:scoremean}
  d_x^{-1}\tau_{U_n(j)}G^n_{T_\varepsilon^M\setminus\{U_n(j)\}}(U_n(j),U_n(j)).
\end{equation}

Since Conditions \ref{c:shallow} and \ref{c:deep} ensure that 
the visits of deep traps determine the behaviour of the time change $S_n(j)$, 
the sum $\sum_{i=1}^{j-1} s_n(i)$ can be considered as a good 
approximation of $S_n(r_n(j))$. The next condition guarantees that 
the scores $s_n(i)$ have a good asymptotic behaviour, i.e. are independent and 
have the right tail. 
\begin{condition}
  \label{c:conv} 
  Let $(s_\infty(i):i\in \mathbb N)$ be an i.i.d.~sequence given by 
  $s_\infty(i):=\hat e_i \sigma_\varepsilon^M(i)$, where 
  $\sigma_\varepsilon^M(i)$ is a sequence of 
  i.i.d.~random variables taking values between $\varepsilon $ and $M$ 
  with common distribution function
  \begin{equation}
    \label{e:defsigman}
    \mathbb P[\sigma_\varepsilon^M(i)\le u]=
    \frac {\varepsilon^{-\alpha }-u^{-\alpha }}
    {\varepsilon^{-\alpha }-M^{-\alpha }}=:
    \frac{\varepsilon^{-\alpha }-u^{-\alpha }} {p_\varepsilon^M }
    , \qquad u\in [\varepsilon ,M].
  \end{equation}
  and $\hat e_i$ is an i.i.d.~sequence of exponential, mean-one 
  random variables independent of $\sigma_\varepsilon^M$. Then there 
  exists a constant $\mathcal K>0$ such that for all $\varepsilon$, 
  $M$ and for a.e.~$\boldsymbol \tau $,  the sequence 
  $(s_n(j)/t_w(n),j\in \mathbb N)$ converges as $n\to \infty$ in law to 
  the sequence $(\mathcal K s_\infty(j),j\in \mathbb N)$. (For 
    notational convenience we define $s_n(j)=s_\infty(j)$ for all 
    $j\ge \zeta_n$.)
\end{condition}

The last three conditions will ensure that the approximation by the 
$\alpha $-stable 
subordinator is relevant for the aging. 

First, we need that $S_n(r_n(\zeta_n))$ is 
larger than $(1+\theta)t_w(n)$. Since 
$r_n(\zeta_n)\ge \sum_{i=1}^{\zeta_n-1}s_n(i)$, and $s_n(i)$ are 
easier to control than $S_n(r_n(j))$ we require
\begin{condition}
  \label{c:large}
  For a.e.~$\boldsymbol \tau $ and for any fixed 
  $\theta >0$, $\delta >0$ it is possible to choose $\xi_n$ such that for all 
  $\varepsilon $ small and $M$ large enough, and for $\zeta_n$ 
  defined in \eqref{e:defzeta}
  \begin{equation}
    \label{e:c:large}
    \mathbb P\Big[\sum_{i=1}^{\zeta_n-1}s_n(i)\ge (1+\theta 
	)t_w(n)\Big|\boldsymbol \tau \Big]\ge 1-\delta .
  \end{equation}
\end{condition}

Second, to prove aging for the two-point function $R$  we need to show that for 
any time $t'$ between $S_n(r_n(j))$ and $S_n(r_n(j+1))$ the probability that 
$X_n(t')=U_n(j)$ is large. 
For a formal statement of this claim we need  some definitions.
Let $t'_n$ be a deterministic time sequence satisfying 
$t_w(n)/2\le t'_n\le (1+\theta )t_w(n)$, and let $\delta >0$.
We define $j_n\in \mathbb N$ by 
\begin{equation}
  \label{e:barjn}
  S_n(r_n(j_n))\le t'_n\le S_n(r_n(j_n+1))-\delta t_w(n), 
\end{equation}
and $j_n=\infty$ if \eqref{e:barjn} is not satisfied for any integer.
Let $A_n(\delta )$ be an event defined by
\begin{equation}
  \label{e:defAn}
  A_n(\delta )=\{ 0< j_n<\zeta_n\}.
\end{equation}
We require
\begin{condition}
  \label{c:post}
  For any $\delta$ it is possible to choose $\varepsilon $ small  and $M$ large 
  enough such that for a.e.~$\boldsymbol \tau $ and all $n$ large enough
  \begin{equation}
  \mathbb P[X_n(t'_n)=U_n(j_n)|A_n(\delta ),\boldsymbol \tau ]\ge 1-\delta.
  \end{equation}
\end{condition}

The last condition that we need excludes repetitions in the sequence $U_n$.
\begin{condition}
  \label{c:noreturn}
  \label{c:last}
  For any fixed $\varepsilon$ and $M$ and a.e.~$\boldsymbol \tau $
  \begin{equation}
    \lim_{n\to \infty}\mathbb P[\exists 0< i,j\le \zeta_n \text{ such 
	that } i\neq j \text { and } U_n(i)=U_n(j)|\boldsymbol\tau  ]=0.
  \end{equation}
\end{condition}

\medskip

We have formulated the six conditions that are inspired by the 
complete graph proof. It should be then not surprising that they 
imply the same result as for the complete graph:

\begin{theorem}[Aging on general graphs]
  \label{t:main}
  Assume that Conditions \ref{c:shallow}--\ref{c:last} holds.
   Then for a.e.~realisation of the random environment 
  $\boldsymbol \tau $
  \begin{equation}
    \lim_{n\to\infty}R_n(t_w(n),(1+\theta )t_w(n);\boldsymbol \tau 
      )=\Asl_\alpha (1/1+\theta ).
  \end{equation}
\end{theorem}

This theorem can be proved in a very similar way as 
Theorem~\ref{t:agingcg}. The complete proof can be found in 
\cite{BC06}. 

\begin{remark}
  1. The set $T_\varepsilon^M(n)$ as defined in \eqref{e:deftop} is 
  sometimes too large and not all conditions that we formulated can be 
  verified easily. Typically this happens when  
  two points of the top are too close to each other with  non-negligible probability. In 
  such case  it is useful to define a set 
  $B(n)\subset T_\varepsilon^M(n)$ of bad traps which cause the 
  difficulties. If it is then possible to verify 
  Conditions~\ref{c:shallow}--\ref{c:last} for the set 
  ${T'}_\varepsilon^M(n):=T_\varepsilon^M(n)\setminus B(n)$ and, 
  moreover, if $B(n)$ satisfies for a.e. $\boldsymbol \tau $ a similar 
  condition as the set of the very deep traps $T_M$, 
  \begin{equation}
    \limsup_{n\to \infty} \mathbb P[H_n(B(n))\le \xi_n|\boldsymbol \tau ]= 0,
  \end{equation}
  then the conclusions of Theorem~\ref{t:main} hold without change.

  2. The last pair of conditions is, in principal, necessary only for a 
  ``post-processing''. If they are not verified, it is possible to  prove
  aging for a {\em top-dependent} correlation function
  \begin{equation}
    \label{e:Rprime}
    R'_n(t_w,t_w+t;\boldsymbol \tau )
    =\mathbb P\big[\exists j :
      S_n(r_n(j))\le t_w < t_w +t\le S_n(r_n(j+1))\big|\boldsymbol \tau \big],
  \end{equation}
  which gives the probability that at most one site in the top is 
  visited by $X_n$ during the observed time interval. A two-point 
  function similar to $R'_n$ was considered in \cite{BBG03}.
\end{remark}

\begin{theorem}
  If only Conditions \ref{c:shallow}--\ref{c:large} hold, then 
  $\boldsymbol \tau $-a.s.
  \begin{equation}
    \lim_{n\to\infty}R'_n(t_w(n),(1+\theta )t_w(n);\boldsymbol \tau 
      )=\Asl_\alpha (1/1+\theta ).
  \end{equation}
\end{theorem}

\subsection{Potential-theoretic characterisation}
\label{ss:verif}
In the previous paragraph we have stated six conditions that allow 
to prove aging on an arbitrary sequence of graphs. It is however not 
clear if these conditions can be verified for any concrete model. In 
fact, they are satisfied for the graphs: $G_n=G=\mathbb Z^d$ with 
$d>1$; for $G_n$ a large torus in $d$ dimensions, 
$G_n=\mathbb Z^d/n\mathbb Z^d$; for the $n$-dimensional hypercube, 
$G_n=\{-1,1\}^n$ among others (included of course the complete 
  graph!). These examples will be developed in the next section, but 
we wont be able to give proofs, which can be pretty difficult (see 
  \cite{BCM06, BC06}). Rather, we want to give here a few hints with 
potential-theoretic flavor on how to verify our conditions.

We want mainly discuss the crucial Condition~\ref{c:conv}: the 
convergence of the scores $s_n(i)$ to the i.i.d.~sequence $s_\infty$. 
In the discussion we suppose that $G_n$ are finite and sufficiently 
regular. Observe first that the set of deep traps $T_\varepsilon^M$ 
is a \textit{random cloud} on $\mathcal V_n$ (i.e. set of points 
  chosen independently from $\mathcal V_n$, we assume that $\tau_x$ 
  are i.i.d.). The intensity $\rho_n$ of this cloud depends on $n$;  
under Assumption~\ref{a:limit} 
$\rho_n\sim p_\varepsilon^M g(n)^{-\alpha } = (\varepsilon^{-\alpha }-M^{-\alpha })g(n)^{-\alpha }$. 
We, obviously, need that $\rho_n\to 0$, so that the random cloud 
is sufficiently sparse, but also $\rho_n |\mathcal V_n|\to \infty$, 
so that the mean size of the random cloud diverges (otherwise 
  fluctuations of the random depths  are important and no 
  a.s.~convergence holds). 

Going back to sequence $s_n(i)$, we have already remarked that
conditionally on $U_n(j)$, the $j^{\text{th}}$ visited deep trap, the 
score $s_n(j)$ is exponentially distributed with mean
\begin{equation}
  \label{e:scoremeana}
  \nu_n^{-1}d_{U_n(j)}^{-1}
  \tau_{U_n(j)}G^n_{T_\varepsilon^M(n)\setminus\{U_n(j)\}}(U_n(j),U_n(j)),
\end{equation}
where $G^n_A(\cdot,\cdot)$ is the Green's function of the simple random 
walk $Y_n$ that is killed on the first hit of the set 
$A\subset\mathcal V_n$. 

If  the graphs $G_n$ are ``sufficiently regular'',  then
all vertices have similar degree,  that means, e.g.,~that there is a 
scale $d(n)$ such that $d_x/d(n)$ is uniformly bounded from 0 and 
$\infty$. Then by setting $\nu_n=d(n)^{-1}$ we can, at least 
theoretically, ignore the first two terms in \eqref{e:scoremeana}. This 
problem does not appear in all examples we consider, there always 
$d_x=\text{const}(n)$. Hence, it remains to control $\tau_{U_n(j)}$ 
and the Green's function in \eqref{e:scoremeana}.  

For the Green's function, one typically prove that there is a scale 
$f(n)\sim t_w(n)/g(n)$ such that $\boldsymbol \tau $-a.s. 
\begin{equation}
  \label{e:greenhom}
  f(n)^{-1}G^n_{T_\varepsilon^M(n)\setminus\{x\}}(x,x)
  \xrightarrow{n\to\infty}\text{const.}
\end{equation}
uniformly for all $x\in T_\varepsilon^M(n)$. This result is again 
reasonable if the graph $G_n$ is sufficiently regular and finite, and 
the cloud $T_\varepsilon^M(n)$ is very diluted. 

To control the distribution of $\tau_{U_n(j)}$ consider first an 
arbitrary random cloud $A_n\subset\mathcal V_n$ with intensity 
$c \rho_n$. This random cloud will represent the set of the deep 
traps or its subsets.  Recall that $H_n(A)$ denotes the hitting time 
of $A$ by $Y_n$. Let $\mathbb P_x$ be the law of $Y_n$ started at 
$x$. Suppose that it is possible to show for all $u>0$ and some scale 
$r(n)$ independent of $c$
\begin{equation}
  \label{e:unifexp}
  \sup_{x\in A_n}
  \Big|\mathbb P_x\Big[
    \frac{ H_n(A_n\setminus\{x\})}{r(n)}\ge u \Big] 
  - \exp(- c u)\Big|
  \xrightarrow{n\to\infty}0, \qquad \boldsymbol \tau -\text{a.s.},
\end{equation}
so that the distribution of the normalised hitting time converges 
uniformly to the exponential distribution with mean $c^{-1}$. This is 
again a reasonable  property for very diluted clouds.  If 
\eqref{e:unifexp} holds for all $c$ then the lack-of-memory property 
of the exponential distribution allows to prove the following claim: 
Let $A_n$ be a random cloud with intensity $(a+b)\rho_n$ and let 
$B_n\subset A_n$ be its sub-cloud with intensity $b \rho_n$. Then 
uniformly for all $x\in A_n$ 
\begin{equation}
  \mathbb P_x[H(B_n\setminus\{x\})<H(A_n\setminus\{x\})]
  \xrightarrow{n\to\infty} \frac b{a+b}.
\end{equation}
This claim yields, e.g.~under Assumption~\ref{a:limit}, that for 
$u\in (\varepsilon ,M)$
\begin{equation}
  \mathbb P_{U_n(j-1)}[\tau_{U_n(j)}\le u]
  \xrightarrow{n\to\infty} 
  \frac 
  {p_\varepsilon^u}{p_\varepsilon^M},
\end{equation}
which is exactly what gives the ``right tail'' and the asymptotic 
independence of $s_n(i)$. The Condition~\ref{c:conv} is then 
consequence of the claims of the last three paragraphs.

If the graphs $G_n$ are infinite one typically cannot prove 
uniformity in \eqref{e:greenhom} and \eqref{e:unifexp}.  One can 
however restrict only to those deep traps that are reachable in 
$\xi_n$ steps from the starting position. On this set the uniform 
control is usually possible.  

There is also a heuristic reason for Condition~\ref{c:conv}. As 
the random cloud $T_\varepsilon^M(n)$ becomes more diluted, the hitting
measure of this cloud charges more and more points because the random walk 
can ``pass more easily around''. If this happens, the 
law-of-large-numbers-type arguments hold not only for the whole sets
$T_\varepsilon^u(n)$ but also for these sets \textit{as sampled by the 
  hitting measure}. That is why Condition~\ref{c:conv} holds.

\begin{remark}
  This is \textit{not true} on $\mathbb Z$! The simple random walk on 
  $\mathbb Z$ cannot ``pass around'': it hits necessarily one of two 
  points of the random cloud that are ``neighbours'' of its starting 
  position. This explains the special properties of the BTM on 
  $\mathbb Z$ and also the necessity of averaging in 
  Theorems~\ref{TF:AGING} and \ref{tf:subaging}. 
\end{remark}

The Conditions~\ref{c:deep} and~\ref{c:large} are also easy 
consequences of claim~\eqref{e:unifexp}. One need first fix 
$\xi_n = m r(n)$ with $m$ large enough to accumulate large enough 
(but finite) number $K$ of ``independent'' scores $s_n(i)$ in order to have 
with a large probability $\sum_{i=1}^K s_n(i)>1+\theta $. This 
satisfies Condition~\ref{c:large}. Then one fix $M$ large enough, 
such that the intensity  $c(M) \rho_n$ of $T_M(n)$ satisfies 
$c(M)>\delta^{-1} m$. Then Condition~\ref{c:deep} holds.

The preceding discussion can be summarised as follows:
\begin{claim}
  If \eqref{e:greenhom} and 
  \eqref{e:unifexp} can be checked for $T_\varepsilon^M(n)$ on $G_n$, then, 
  under Assumption~\ref{a:limit},  
  Conditions~\ref{c:deep}--\ref{c:large} can be verified.
\end{claim}

The Condition~\ref{c:shallow} does not follow directly from 
\eqref{e:unifexp} since Poisson clouds with larger intensity than 
$\rho_n$ should be considered. The ``slicing strategy'' as presented 
in the proof for the complete graph however usually works. The 
remaining Conditions~\ref{c:post} and~\ref{c:noreturn} are not 
substantial and we do not discuss them here.

Remark also that Condition~\ref{c:shallow} together 
with~\eqref{e:greenhom} and~\eqref{e:unifexp} allows to prove the 
approximation of the clock process by an $\alpha $-stable 
subordinator. This approximation is not a consequence of Conditions 
1--6 only.

There are at least two methods of  proof for  
facts~\eqref{e:greenhom} and~\eqref{e:unifexp}. The first is the 
coarse-graining procedure of \cite{BCM06} that is explained in 
Section~\ref{ss:coarse}. The advantage of this procedure is that it 
should work on many different graphs. It is however relatively 
technical and many special cases should be treated apart.  One can 
also use the formula given by Matthews \cite{Mat88}. This method is 
used in \cite{BC06} to prove aging in the REM and on the torus in 
$\mathbb Z^2$. It however applies only if $Y_n$ ``completely 
forgets'' the position of $U_n(i)$ before hitting $U_n(i+1)$, 
i.e.~that the hitting measure of the random cloud is essentially 
uniform. Therefore this method does not apply e.g.~on $\mathbb Z^d$.

\section{Applications of the arcsine law}
\label{s:applications}
We now describe two examples where the approach of 
Section~\ref{s:general} can be used to prove aging.

\subsection{Aging in the REM} 
The Random Energy Model is  the 
simplest mean-field model for  spin-glasses and its static behaviour 
is well understood. The studies of dynamics are much more sparse. The 
first proof of aging in the REM was given in  \cite{BBG03,BBG03b}, 
based on 
renewal theory. The general approach of 
Section~\ref{s:general} gives another, shorter proof. This approach 
allows to prove aging on a broader range of time scales, but on the 
other hand do not include quite exactly the results of 
\cite{BBG03,BBG03b}. We will compare both results later. Before doing 
it, let us define the model and give some motivation why and in which 
ranges of times and temperatures aging occurs. 

The Random Energy model is a mean-field model of a spin-glass. It 
consists of $n$ spins that can take values $-1$ or $1$, that is 
configurations of the REM are elements of $\mathcal V_n=\{-1,1\}^n$. The 
energies $\{E_x,x\in \mathcal V_n\}$ of the configurations are 
i.i.d.~random variables. The standard choice of the marginal 
distribution of $E_x$ is centred normal distribution with variance 
$n$. We will, however, deviate from  the standard choice to simplify 
the computations and we will assume that $-E_x$ are i.i.d.~positive 
random variables with the common distribution given by
\begin{equation}
  \label{e:Edist}
  \mathbb P[-E_x/\sqrt n\ge u]= e^{-u^2/2},\qquad u\ge 0.
\end{equation}
We then define
\begin{equation}
  \tau_x=\exp(-\beta  E_x).
\end{equation}
The distribution \eqref{e:Edist} has almost the same tail behaviour as 
the normal distribution. As we already know, it is the tail behaviour 
of $\tau_x$ (and thus of $E_x$)  that is responsible for aging. 
Therefore, the 
use of ``faked normal distribution'' \eqref{e:Edist} is not
substantial for our discussion. Remark also that a similar trick, 
i.e. to take $E_x$ to be minimum of $0$ and the normal variable, was 
used in \cite{BBG03b}. 

For the dynamics of the REM we require that only one spin can be 
flipped at a given moment. This corresponds to 
\begin{equation}
  \mathcal E_n = \{\<x,y\>\in \mathcal V_n^2:\sum_{i=1}^n|x_i-y_i|=2\},
\end{equation}
where  $(x_1,\dots,x_n)$ are the values of individual spins. We use $G_n$ 
to denote the $n$-dimensional hypercube $(\mathcal V_n,\mathcal E_n)$. 
There are many choices for the dynamics of REM, that has the Gibbs 
measure $\boldsymbol \tau $ as a reversible measure. We will 
naturally consider the trap model dynamics $\BTM(G_n,\boldsymbol \tau ,0)$, 
which is one of the simplest choices. We fix $\nu_n=1/n$ in 
definition \eqref{e:rates}, so that $\tau_x$ is the mean waiting time at $x$.
We always suppose that
\begin{equation}
  Y_n(0)=X_n(0)=\boldsymbol 1=(1,\dots,1).  
\end{equation}

\subsubsection{Short time scales}
\begin{theorem}
  \label{t:hcaging}
  Let the parameters $\alpha \in (0,1)$ and $\beta >0$ be such that 
  \begin{equation}
    \label{e:hcagass}
    3/4 < \alpha^2 \beta^2 / 2\log2 < 1.
  \end{equation}
  Define
  \begin{equation}
    t_w(n):= \exp(\alpha \beta^2 n).
  \end{equation}
  Then, for a.e.~$\boldsymbol \tau $,
  \begin{equation}
    \lim_{n\to \infty}
    R_n(t_w(n),(1+\theta )t_w(n))=
    \Asl_\alpha (1/1+\theta ).
  \end{equation}
\end{theorem}

Let us first explain the appearance of scale $t_w(n)$ together with 
one problem that is specially related to REM.
We have seen in 
Section \ref{s:general} that aging occurs only if $\tau_x$ are 
sufficiently heavy-tailed. This certainly fails to be true 
for the REM: an easy calculation gives 
$\mathbb P[\tau_x\ge u]=u^{-\log u/2\beta^2 n}$, which decreases 
faster than any polynomial. It is therefore clear that, if the system 
is given enough time to explore a large part of the configuration 
space and thus to discover the absence of heavy tails, then no aging 
occurs, at least not in our picture. On the other hand, at shorter time 
scales the system does not feel the non-existence of heavy tails as 
can be seen from the following estimate. Let $\alpha >0$, then
\begin{equation}
  \begin{aligned}
    e^{\alpha^{2}\beta^2 n/2} 
    \mathbb P\big[\tau_x\ge  u e^{\alpha \beta^2 n}\big]
    &=
    e^{\alpha^{2}\beta^2 n/2} \mathbb 
    P\Big[E_x\ge \frac {\log u + \alpha \beta^2 n}
      {\beta \sqrt{n}}\Big]\\
    &=
    \exp\Big\{-\frac{\log^2 u}{2\beta^2 n}-\alpha \log u\Big\}
    \xrightarrow{n\to\infty} u^{-\alpha }.
  \end{aligned}
\end{equation}
and therefore
\begin{equation}
  \label{e:htcomp}
  \mathbb P\Big[\frac {\tau_x}{e^{\alpha \beta^2 n}}\ge u\Big]=
    e^{-\alpha^2\beta^2n/2}\cdot u^{-\alpha }(1+o(1))\qquad (n\to\infty).
\end{equation}
In view of the fact that the simple random walk on the hypercube 
almost never backtracks, it seems reasonable to presume that 
if the process had time to make only approximately 
$e^{\alpha^2 \beta^2 n/2}$ steps,  then it has no time to discover 
the absence of heavy tails and aging could be observed. The above 
theorem shows this presumption to be true.

Let us remark that there is much stronger relation between ``random 
exponentials'' $\tau_x$ and heavy-tailed random variables. Let 
$(F_i,i\in \mathbb N)$ be an i.i.d.~sequence of centred normal 
random variables with  variance one.  It was proved in 
\cite{BBM05} that for some properly chosen $Z(n)$ and $N(n)$ the 
normalised sum 
\begin{equation}
  \frac 1 {Z(n)} \sum_{i=1}^{N(n)}e^{-\beta \sqrt n F_i}
\end{equation}
converges as $n\to \infty$ in law to an $\alpha $-stable distribution 
with $\alpha $ depending on $\beta$ and  $N(n)$. We show that  the 
same is true for the properly normalised clock process $S(n)$, 
which is a properly normalised sum of correlated random variables, 
more precisely of i.i.d.~random variables sampled by a random walk.

In view of \eqref{e:htcomp} it is easy to fix  objects for which  
Conditions~\ref{c:shallow}--\ref{c:last} should be verified:  
we define
\begin{align}
  \label{e:hcparmsa}
  t_w(n)&:=
  \exp(\alpha \beta^2 n),\\
  \label{e:hcparmsb}
  \xi_n& := m \exp(\alpha^2 \beta^2 n/2), \\
  T_\varepsilon^M(n,\alpha )&:=\{x\in \mathcal V_n: 
    \tau_x\in(\varepsilon ,M)e^{\alpha \beta^2 n}\}. 
  \label{e:hcparmsc}
\end{align}

The Theorem~\ref{t:hcaging} is then the consequence of the following 
proposition and Theorem~\ref{t:main}.
\begin{proposition}
  \label{p:hcaging}
  Let $\alpha $ and $\beta $ be as in Theorem~\ref{t:hcaging}. Then 
  for any $\theta $ it is possible to choose $m$ large enough such 
  that Conditions~\ref{c:shallow}--\ref{c:last} hold for 
  $\mathbb P$-a.e.~$\boldsymbol \tau $. 
\end{proposition}
We believe that the range of the validity \eqref{e:hcagass} of the 
Theorem~\ref{t:hcaging} is not the broadest possible.  The upper 
bound $1$ is correct. If $\alpha^2 \beta^2 / 2\log2>1$, then 
$\xi_n\gg 2^n$.  That means that the state space $\mathcal V_n$ 
becomes too small and the process can feel its finiteness.   
On the other hand, the lower-bound $3/4$ is 
purely technical and can probably be improved. 

Observe also that the condition \eqref{e:hcagass}  can be rewritten as
\begin{equation}
  \alpha^{-1}\beta_c \sqrt{3/4} < \beta  < \alpha^{-1}\beta_c,
\end{equation}
where $\beta_c = \sqrt{2\log 2}$ is the critical temperature of the 
usual REM. This, in particular, means that aging can be observed in 
REM also \textit{above the critical temperature}, $\beta < \beta_c$.  

The proof of this proposition in \cite{BC06} follows the strategy 
outlined in Section~\ref{s:general} and uses
the results of Matthews \cite{Mat88} for the fine 
control of the simple random walk on the hypercube. In particular, 
\eqref{e:unifexp} is a consequence of the following 
potential-theoretic result which might be of independent interest.
\begin{proposition}
  \label{p:hcpot}
  (i) Let for all $n\ge 1$ sets $A_n\subset \mathcal V_n$ be such that 
  $|A_n|=\rho_n 2^n$ with ``densities'' $\rho_n$ satisfying 
  $\lim_{n\to \infty}\rho_n2^{\gamma n}=\rho \in (0,\infty)$ for some
  $\gamma \in (1/2,1)$. Let further the sets $A_n$ satisfy the {\em minimal 
    distance condition}
  \begin{equation}
    \label{e:distcond}
    \min\{d(x,y):x,y\in A_n\}\ge (\omega(\gamma )  +\varepsilon )n
  \end{equation}
  for some small constant $\varepsilon>0$ and for the unique solution 
  $\omega(\gamma )$  of 
  \begin{equation}
    \omega\log \omega+(1-\omega)\log(1-\omega) +\log 2=(2\gamma -1)\log 2,
    \qquad \omega \in(0,1/2).
  \end{equation}
  Then for all $s\ge 0$
  \begin{equation}
    \lim_{n\to\infty}\max_{x\in A_n}\Big |
    \mathbb E_x \Big[\exp\Big(-\frac s {2^{\gamma 
	    n}}H_n\big(A_n\setminus\{x\}\big)\Big)\Big]-
    \frac \rho {s+\rho }\,\Big|=0.
  \end{equation}
  That means that the hitting time $ H_n(A_n\setminus\{x\})/2^{\gamma n}$ is asymptotically 
  exponentially distributed  with mean~$1/\rho $.

  (ii) If $A_n$ are random clouds with intensity $\rho_n$ such that 
  $\lim_{n\to \infty}\rho_n2^{\gamma n}=\rho \in (0,\infty)$ and 
  $\gamma \in (3/4,1)$, then the assumptions of (i) are a.s.~satisfied. 
\end{proposition}

This result also explains the appearance of the lower bound $3/4$ in 
the range of the validity of Theorem~\ref{t:hcaging}: the set of deep 
traps satisfies the assumptions of Proposition~\ref{p:hcpot}(ii) only 
if $\alpha^2 \beta^2/2\log 2>3/4$.

\subsubsection{Long time scales}
The result of \cite{BBG03,BBG03} deals with the longest possible 
time scales where aging appears in REM. 
 The continuous-time 
Markov process $X$ is replaced by a discrete-time process $X'$, which 
at every step has the possibility not to move. The number of tries before 
leaving $x$ has geometrical distribution with mean $\tau_x$.  
As $n\to \infty$, this dynamics differs very little from the usual 
trap model dynamics.
The random environment is given by
\begin{equation}
  \tau_x=\exp(\beta\sqrt n \max(E_x,0)),
\end{equation}
where $E_x$ are i.i.d. centred normal random variables with variance 
one. 

To define the relevant time and depth scales we set
\begin{equation}
  u_n(E)=\beta_c \sqrt n + E/\beta_c\sqrt n - \log( 4\pi n\log 2) /2\beta_c \sqrt n.
\end{equation}
We define the set of deep traps (\textit{the top}) by
\begin{equation}
  T_n(E)=\{x\in \mathcal V_n: E_x\ge u_n(E)\}.
\end{equation}
The function $u_n(E)$ is chosen in such way that, as $n\to\infty$ and 
$E$ is kept fixed the distribution of $|T_n(E)|$ converges to the 
Poisson distribution with mean that depends only on  $E$. 
The mean diverges if $E\to -\infty$ afterwards. The initial position 
of the discrete-time Markov chain $X'$ is chosen to be uniformly 
distributed in $T_n(E)$. 

A different correlation function considered in \cite{BBG03}:   
If $x_n(k)$ denotes the last trap visited by $X'$ before step $k$ 
then 
\begin{multline}
  \Pi'_n(k,k+l,E;\boldsymbol \tau )
  = 
  \mathbb P\big[\{X'(i):i\in\{k+1,\dots, k+l\}\cap (T_n(E)\setminus\{x_n(k)\})=
      \emptyset\big|\boldsymbol \tau \big].
\end{multline}
It is essentially the same as the function $R'$ (see \eqref{e:Rprime})  

This is the main aging result of \cite{BBG03b}
\begin{theorem}
  \label{t:bbgaging}
  For any $\beta > \beta_c=\sqrt{2\log 2}$ and any $\varepsilon >0$
  \begin{equation}
    \label{e:BBGresult}
    \lim_{t\to\infty}\lim_{E\to -\infty}\lim_{n\to \infty}
    \mathbb P\Big[\Big|\frac{\Pi'_n(c_n t, (1+\theta )c_n t, 
	  E;\boldsymbol \tau )}{\Asl_{\beta_c/\beta }(1/1+\theta )} -1 \Big|> \varepsilon \Big]=0,
  \end{equation}
  where $c_n\sim e^{\beta \sqrt n u_n(E)}$.
\end{theorem}
 
In fact one can see that the dynamics of the REM when observed only 
on the top $T_n(E)$ can be approximated very well when $n\to\infty$ 
and $E\to -\infty$ by a BTM on the complete graph with 
$M=|T_n(E)|\sim e^{-E}$ vertices (see \cite{BBG02,BBG03}).

\medskip

Let us now compare the results of Theorems~\ref{t:hcaging} 
and~\ref{t:bbgaging}. First, 
different correlation functions $R$ and $\Pi '$  are considered. 
This difference is not substantial, we believe that it 
is possible to eliminate the top dependence (i.e.~to convert 
  something like $R'$ to something like $R$) of \eqref{e:BBGresult} 
by some post-processing in the direction of Condition~\ref{c:post}. 

The a.s.~convergence in Theorem~\ref{t:hcaging} is stronger than the 
convergence in probability in Theorem~\ref{t:bbgaging}. It is a 
consequence of the fact that much larger time scales are considered 
and the set of the deep traps $T_n(E)$ is finite for fixed $E$. 
Therefore, we cannot use law-of-large-numbers-type arguments for the 
number of deep traps with depth in a fixed interval. We have seen 
this effects already in the case of the complete graph (see 
  \eqref{e:cginproba}).

The main difference between the two theorems  is in the considered 
top sizes and time scales. In Theorem~\ref{t:bbgaging} the size of 
the  top is kept bounded as $n\to \infty$. This allows to apply 
``lumping techniques'' to describe the properties of the projection 
of a simple random walk on the hypercube to the top, that is to prove 
that that an equivalent of the process $U_n$ (see \eqref{e:Undef}) 
converges to the simple random walk on the complete graph with the 
vertex set $T_n(E)$.  In Theorem~\ref{t:hcaging} (or more precisely,  
  in Proposition~\ref{p:hcaging}) the size of the top 
$T_\varepsilon^M(n)$ increases exponentially with $n$. This makes the 
application of the lumping not possible and techniques based on 
Matthews results, i.e.~Proposition~\ref{p:hcpot}, should be used. 

The time scale 
$c_n\sim e^{\beta \sqrt n u_n(E)}\sim e^{\beta \beta_c n + \beta E/\beta_c}$ 
of Theorem~\ref{t:bbgaging} corresponds to the case 
$\alpha \beta/\beta_c =1$ and is much larger than the scale 
$t_w(n)\sim e^{\alpha \beta^2n}=e^{\frac{\alpha \beta }{\beta_c}\beta \beta_c n }$. 
These scales approach if $\alpha \beta /\beta_c$ tends to $1$, which 
is the upper limit of the validity of Theorem~\ref{t:hcaging}. It 
would be possible to improve this theorem by setting
$t_w(n)=e^{\beta \beta_c n} f(n)$ with some $f(n)\to 0$ as 
$n\to \infty$ sufficiently fast, but even then $t_w(n)\ll c_n$. 
Exactly at $\alpha \beta /\beta_c=1$ Theorem~\ref{t:hcaging} does not 
hold. As we have already remarked, in this case it is necessary to 
use the double-limit procedure and the convergence in probability.

\subsubsection{Open questions and conjectures}

The Theorems~\ref{t:hcaging} and~\ref{t:bbgaging} give rigorous proofs of aging in the REM. They 
are however only partly satisfactory. It would be nice to replace the 
RHT dynamics by a more physical dynamics, like e.g. Glauber, or, at 
least, to explore the $a\neq 0$ case. We believe that the long-time 
behaviour of the model should not change dramatically, however we 
do not know any proof of it. The problem is that the Markov chain $X$ 
becomes a time change of a random walk in random environment on the 
hypercube. Moreover, the clock process and the random walk are 
dependent. 

Another natural direction of research is to extend the results for 
the RHT dynamics on the REM to other mean-field spin-glasses, like 
the SK model or the $p$-spin SK model. In these models the energies 
of the spin configurations $E_x$ are no longer independent. We 
strongly believe that the approach of Section~\ref{s:general} can be 
applied, at least for  large $p$ and for well chosen time scales. 
These scales should be short enough not to feel the extreme values of 
the $E_x$'s which rule the (model-dependent) statics, but long enough 
for the convergence to  a stable subordinator to take place for the 
clock-process. The difficulty is to verify Conditions~1--6 if the 
$E_x$'s are not i.i.d. The assumption that the $E_x$'s are 
independent is used twice in the proof for the REM. First, we use it 
to verify Condition~\ref{c:shallow}, that is to prove that the time 
spent in the shallow traps is small. We believe that this condition 
stays valid  also for dependent spin-glass models.  The second use of 
the independence is  more substantial. It is used to describe the 
geometrical structure of the set of the deep traps. More exactly it 
is used  to bound from below the minimal distance between deep traps 
and to show that the number of traps in $T_u^{u+\delta }(n)$ is 
proportional to $u^{-\alpha }$. It is an open question if these 
properties remain valid for dependent spin-glasses.

\subsection{Aging on large tori}
Another graph where the approach of Section~\ref{s:general} can be 
used to prove aging
is a torus in $\mathbb Z^d$. For convenience we will consider only 
$d=2$ here, although similar results are expected to hold for 
$d\ge 3$. 

Let $G_n=(\mathcal V_n, \mathcal E_n)$ be the two-dimensional torus of size 
$2^n$ with nearest-neighbours edges, i.e. 
$\mathcal V_n=\mathbb Z^2/2^n\mathbb Z^2$, and edge $\<x,y\>$ is in 
$\mathcal E_n$ iff 
\begin{equation}
  \sum_{i=1}^2 |x_i-y_i| \bmod 2^n = 1. 
\end{equation}
We 
use $d(x,y)$ to denote the graph distance of $x,y\in \mathcal V_n$. Let 
further $\boldsymbol \tau =\{\tau^n_x\}$, $x\in \mathcal V_n$, 
$n\in \mathbb N$, be a collection of positive i.i.d.~random variables 
satisfying Assumption~\ref{a:limit}.
We consider the Bouchaud trap model, $\BTM(G_n,\boldsymbol \tau ,0)$ with 
$\nu =1/4$.

\begin{theorem}[Aging on the torus]
  \label{t:Zdaging}
  Let $t_w(n)=2^{2n/\alpha }n^{1-(\gamma /\alpha)}$ with $\gamma \in (0,1/6)$.
  Then for $\mathbb P$-a.e.~realisation of the random environment 
  $\boldsymbol \tau $
  \begin{equation}
    \lim_{n\to\infty} R(t_w(n),(1+\theta )t_w(n);\boldsymbol \tau ) 
    =\Asl_\alpha (1/1+\theta ).
  \end{equation}
\end{theorem}

The theorem follows from the following proposition whose proof can be 
found in  \cite{BC06}. 
\begin{proposition}
  \label{p:Zdcond}
  For any $\theta $ there exist $m$ large enough such that Conditions 
  1--6 hold for 
  \begin{equation}
    \label{e:parmZ}
    \begin{gathered}
      t_w(n)=2^{2n/\alpha }n^{1-\gamma /\alpha },\qquad\qquad
      \xi_n =m 2^{2n}n^{1-\gamma },\\
      T_\varepsilon^M(n)=
      \{x\in \mathcal V_n:\tau_x\in (\varepsilon ,M) 
	2^{2n/\alpha }n^{-\gamma /\alpha }\},
    \end{gathered}
  \end{equation}
\end{proposition}

The main motivation for Theorem~\ref{t:Zdaging} was to extend the 
range of aging scales on $\mathbb Z^2$ and mainly to really explore 
the extreme values of the random landscape. Namely, the BTM on the whole 
lattice $\mathbb Z^2$ does not find the deepest traps that are close to its 
starting position. In the first $2^{2n}$ steps, it gets to the distance 
$2^{n}$ and visits $O(2^{2n}/\log(2^{2n}))=O(2^{2n}/n)$ sites. 
Therefore, the deepest visited trap has a depth of order 
$2^{2n/\alpha }/n^{1/\alpha }$, which is much smaller that the depth 
of the deepest trap in the disk with radius $2^n$, that is 
$2^{2n/\alpha }$. Eventually, the process visits also this deepest 
trap, however it will be too late. This trap will no longer be 
relevant for the time change since much deeper traps will have 
already be 
visited. The deepest trap is relevant only if the random walk stays 
in the neighbourhood of its starting point a long enough time. One 
way to force it to stay is to change $\mathbb Z^2$ to the torus. By 
changing the size of the torus relatively to the number of considered 
steps, i.e.~by changing $\gamma $, different depth scales become 
relevant for aging.  

The range of possible values of $\gamma \in (0,1/6)$ has, as in the 
REM case, a natural bound and an artificial one. It is natural that 
$\gamma <0$ cannot be considered, since if the simple walk makes more 
than $2^{2n}(\log2^n)^{1+\varepsilon }$ steps, $\varepsilon >0$, 
inside the torus of size $2^n$, its occupation probabilities are very 
close to the uniform measure on the torus, that is the process is 
almost in equilibrium. The other bound, $\gamma =1/6$, comes from the 
techniques that we use. We do not believe it to be meaningful since we 
expect the theorem to hold for all $\gamma >0$. Actually, the result 
for  $\gamma >1$ follows easily from the Theorem~\ref{t:aging} for 
the whole lattice.  In this case the size of the torus is much larger 
than $\xi_n^2$. So that, the process has no time  to discover 
the difference between the torus and $\mathbb Z^2$. We also know that 
Theorem~\ref{t:Zdaging} holds also in the window $[1/6,1]$ since 
it can be proved by the same methods as for  $\mathbb Z^2$, 
\cite{BCM06}. Nevertheless the complete proof in this window has  
never been written. 

The $\gamma =0$ case corresponds to the longest possible time scales. 
We expect that a similar result as Theorem~\ref{t:bbgaging} is valid. 

The proof of Proposition~\ref{p:Zdcond} uses again Matthews' results: 
the following equivalent of Proposition~\ref{p:hcpot} can be proved.
\begin{proposition}
  \label{p:Zdpot}
  (i) Let $A_n\subset \mathcal V_n$ be such that $|A_n|=\rho_n 2^{2n}$ 
  with the density $\rho_n$ satisfying 
  $\lim_{n\to\infty} 2^{2n}n^{-\gamma }\rho_n=\rho$ 
  for some $\gamma \in (0,1/6)$ and $\rho \in (0,\infty)$. 
  Let further $A_n$ satisfy the minimal distance condition
  \begin{equation}
    \min\{d(x,y):x,y\in A_n\}\ge 2^{n}n^{-\kappa },
  \end{equation}
  for some $\kappa >0 $. Then,
  for $\mathcal K=(2\log 2)^{-1}$,
  \begin{equation}
    \lim_{n\to\infty}\max_{x\in A_n}
    \bigg|\mathbb E_x
    \Big[\exp\Big(-\frac s {2^{2n}n^{1-\gamma }}
	H(A_n\setminus\{x\})\Big) \Big]- 
    \frac {\mathcal K\rho }{s+\mathcal K\rho }\bigg|=0.
  \end{equation}

  (ii) If $A_n$ are, in addition,  random clouds with the   
  densities given above, then the minimal distance condition is  
  a.s.~satisfied for all $n$ large.
\end{proposition}


\appendix
\section*{\appendixname. Subordinators}
\addcontentsline{toc}{section}{Appendix. Subordinators}
\refstepcounter{section}
\label{ap:subord}
We use frequently the theory of increasing L\'evy processes in 
these notes. We summarise in this appendix the facts that are 
important for us. For a complete treatment of this theory the 
reader is referred to the beautiful book by Bertoin \cite{Ber96}.

\begin{definition}
  We say that $V$ is a L\'evy process if for every 
  $s,t\ge 0$, the increment $V(t+s)-V(t)$ is independent of the process 
  $(V(u),0\le u\le t)$ and has the same law as $V(s)$. That means in 
  particular $V(0)=0$. 
\end{definition}

We work only with the class of 
increasing L\'evy processes, so called {\em subordinators}. There is a 
classical one-to-one 
correspondence between  subordinators and the set of pairs 
$(\dd,\mu )$, where $\dd\ge 0$ and $\mu $ is a measure on $(0,\infty)$, 
satisfying
\begin{equation}
  \label{e:lmcond}
  \int_0^\infty (1\wedge x) \mu (\d x)<\infty.
\end{equation}
The law of a subordinator is uniquely determined by 
the Laplace transform of $V(t)$,
\begin{equation}
  \mathbb E[e^{-\lambda V(t)}]=e^{-t \Phi (\lambda )},
\end{equation}
where the {\em Laplace exponent}
\begin{equation}
  \Phi (\lambda )=\dd + \int_0^\infty (1-e^{-\lambda x})\mu (\d x).
\end{equation}
The constant $\dd$ corresponds to the deterministic constant drift. 
All processes  appearing in these notes  have no drift, therefore we 
suppose always $\dd\equiv 0$. The measure $\mu$ is called the {\em 
  L\'evy measure} of the subordinator $V$.

There are two important families of subordinators. The first consists 
of the stable subordinators. A subordinator is stable with 
index $\alpha \in (0,1)$ if for some $c>0$ its Laplace exponent satisfies
\begin{equation}
  \Phi (\lambda )=c\lambda^{\alpha }=
  \frac {c\alpha}  {\Gamma (1-\alpha )} 
  \int_0^{\infty} (1-e^{-\lambda x})x^{-1-\alpha }\,\d x.
\end{equation}
Here $\Gamma $ is the usual Gamma-function. 

The second important family of subordinators are the compound Poisson 
processes. They correspond to finite L\'evy measures,  
$\mu ((0,\infty)) <\infty$. In this case $V$ can be constructed from 
a Poisson point process $J$ on $(0,\infty)$ with constant intensity 
$Z_\mu :=\mu ((0,\infty))$ and a family of i.i.d.~random variables 
$s_i$ with marginal $Z^{-1}_\mu \mu $ as follows. Let 
$J=\{x_i,i\in \mathbb N\}$, $x_1<x_2<\dots$, and $x_0=0$. Then $V$ is 
the process with $V(0)=0$ that is constant on all intervals 
$(x_i,x_{i+1})$, and at $x_i$ it jumps by $s_i$, i.e.~
$V(x_i)-V(x_i-)=s_i$.

We need to deduce  convergence of subordinators from the 
convergence of L\'evy measures. 
\begin{lemma}
  \label{l:subconv}
  Let $V_n$ be subordinators with L\'evy measures $\mu_n$. Suppose that 
  the sequence 
  $\mu_n$ converges weakly to some measure $\mu $ satisfying 
  \eqref{e:lmcond}. Then $V_n$ converge to $V$ weakly in the Skorokhod 
  topology on $D=D([0,T),\mathbb R)$ for all final instants $T>0$. 
\end{lemma}
\begin{proof}[Sketch of the proof.]
  To check the convergence on the space of cadlag path $D$ endowed 
  with Skorokhod topology, it is necessary check two facts: (a)~the 
  convergence of finite-dimensional distributions,  and 
  (b)~tightness. To check (a) it is sufficient to look at 
  distributions at one fixed time, since $V_n$ have independent, 
  stationary increments. From the weak convergence of $\mu_n$ it 
  follows that for all $\lambda >0$,
  \begin{equation}
    \begin{aligned}
      \mathbb E[e^{-\lambda V_n(t)}]=
      e^{-t\int (1-e^{-\lambda x})\mu_n(\d x)}\xrightarrow{n\to\infty}
      e^{-t\int (1-e^{-\lambda x})\mu(\d x)}=
      \mathbb E[e^{-\lambda V(t)}],
    \end{aligned}
  \end{equation}
  which implies the weak convergence of $V_n(t)$. Since $V_n$ are 
  increasing, to check the tightness it is sufficient to check the 
  tightness of $V_n(T)$, which is equivalent to 
  \begin{equation}
    \lim_{\lambda \to 0}\lim_{n\to \infty} \mathbb E[e^{-\lambda V_n(T)}]=1.
  \end{equation}
  This is easy to verify using the weak convergence of $\mu_n$ and 
  the validity of \eqref{e:lmcond} for~$\mu $. 
\end{proof}

\begin{definition}[The generalised arcsine distributions]
  For any $\alpha \in (0,1)$, the {\em generalised arcsine distribution with 
    parameter $\alpha $} is the 
  distribution on $[0,1]$ with density
  \begin{equation}
    \frac{\sin \alpha \pi }{\pi } u^{\alpha -1}(1-u)^{-\alpha }.
  \end{equation}
  We use $\Asl_\alpha $ to denote its distribution function,
  \begin{equation}
    \Asl_\alpha (u):=
    \int_{0}^{u}
    \frac{\sin \alpha \pi }{\pi } u^{\alpha -1}(1-u)^{-\alpha }\,\d 
    u,\qquad u \in [0,1].
  \end{equation}
\end{definition}
Note that 
$\Asl_\alpha (z)=\pi^{-1}\sin (\alpha\pi ) B(z;\alpha ,1-\alpha )$ 
where $B(z;a,b)$ is the incomplete Beta function. It is easy to see that
\begin{equation}
  \label{e:aslass}
  \Asl_\alpha (z)=1-O((1-z)^{1-\alpha }),\qquad\text{as $z\to 1$}.
\end{equation}

The following fact is crucial for us. 
\begin{proposition}[The arcsine law]
  \label{p:stablejump}
  Let $V$ be an $\alpha $-stable subordinator and let 
  $T(x)=\inf\{t:V(t)>x\}$. Then the random variable $V(T(x)-)/x$ has 
  the generalised  arcsine distribution with parameter $\alpha $.
\end{proposition}

This proposition has an important corollary.
\begin{corollary}
  \label{c:stablejump}
  The probability that the $\alpha $-stable subordinator $V$ jumps 
  over interval $[a,b]$ (i.e. there is no $t\in \mathbb R$ such that 
    $V(t)\in [a,b]$) is equal to
  \begin{equation}
    \mathbb P[V(T(b)-)<a]=
    \Asl_\alpha (a/b).
  \end{equation}
\end{corollary}

\begin{proof}[Sketch of the proof of Proposition~\ref{p:stablejump}.]
  Consider the potential measure $U$ of the subordinator~$V$, 
  \begin{equation}
    U(A)=\mathbb E\Big[\int_0^\infty \bbone\{V(t)\in A\}\,\d t\Big].
  \end{equation}
  Its Laplace transform is given by
  \begin{equation}
    \label{e:lapU}
    \int_0^\infty e^{-\lambda x}U(\d x)
    =\mathbb E\Big[\int_0^\infty e^{-\lambda V(t)} \,\d t\Big]
    =\frac 1{\Phi (\lambda )}.
  \end{equation}
  Define further 
  $\bar \mu (x)=\mu ((x,\infty))$. Then
  \begin{equation}
    \label{e:lapmubar}
    \int_{0}^\infty e^{-\lambda x} \bar \mu (x)\,\d x =
    \Phi (\lambda )/\lambda .
  \end{equation}

  Fix $x>0$. For every $0\le y\le x<z$, we can write
  \begin{equation}
    \mathbb P[V(T(x)-)\in \d y, V(T(x))\in \d z ]
    =U(\d y) \mu (\d z- y).
  \end{equation}
  (For a proof of this intuitively obvious claim see p.~76 of 
    \cite{Ber96}.) Define now $A_t(x)=x^{-1}V(T(tx)-)$ and consider 
  its ``double'' Laplace transform
  \begin{equation}
    \tilde A(q,\lambda )=
    \int_0^\infty e^{-qt}\mathbb E[\exp(-\lambda A_t(x))]\,\d t. 
  \end{equation}
  This Laplace transform can be explicitly calculated. Indeed using 
  \eqref{e:lapU} and \eqref{e:lapmubar} we obtain
  \begin{equation}
    \begin{aligned}
      \tilde A(q,\lambda )
      &=
      \int_0^\infty e^{-qt}\int_0^{tx} 
      e^{-\lambda y/x}\bar \mu (tx- y) U(\d y)\,\d t
      \\&=
      \int_0^\infty \int_0^\infty e^{-\lambda y/x}e^{-q(s+y)/x}
      \bar \mu (s) x^{-1}\, U(\d y) \,\d s
      \\&=
      \frac 1x \cdot \frac 1 {\Phi ((\lambda +q)/x)}\cdot \frac {\Phi (q/x)}{q/x}
      =
      \frac {\Phi (q/x)}{q \Phi ((q+\lambda )/x)}. 
    \end{aligned}
  \end{equation}
  Using that $V$ is $\alpha $-stable, i.e.~$\Phi (x)=cx^{\alpha }$, 
  we get
  \begin{equation}
    \tilde A(q,\lambda )=\frac {q^{\alpha -1}}{(q+\lambda )^\alpha }=
    \int_0^\infty \int_0^t e^{-qt}e^{-\lambda s}
    \frac{s^{\alpha -1}(t-s)^{-\alpha }}
    {\Gamma (\alpha )\Gamma (1-\alpha )} \,\d s\,\d t.
  \end{equation}
  Observing that 
  $\pi^{-1}\sin(\alpha \pi ) =(\Gamma (\alpha )\Gamma (1-\alpha ))^{-1}$ 
  yields the claim of the proposition. 
\end{proof}

\addcontentsline{toc}{section}{\refname}

\def\cprime{$'$}
\providecommand{\bysame}{\leavevmode\hbox to3em{\hrulefill}\thinspace}
\providecommand{\MR}{\relax\ifhmode\unskip\space\fi MR }
\providecommand{\MRhref}[2]{%
  \href{http://www.ams.org/mathscinet-getitem?mr=#1}{#2}
}
\providecommand{\href}[2]{#2}

\end{document}